\input ./style/arxiv-general.cfg
\documentclass[aap,MSNbibl,dvips]{arximspdf}
\makeatletter
   \@ifpackageloaded{graphicx}{}{\usepackage{graphicx}}
\makeatother
\usepackage{mathrsfs,mathbh}

\doi{10.1214/15-AAP1133}
\volume{26}
\issue{3}
\pubyear{2016}
\firstpage{1774}
\lastpage{1806}
\docsubty{FLA}

\makeatletter

\newcommand{\eqref}[1]{(\ref{#1})}
\newtheorem{theorem}{Theorem}[section]

\newtheorem{lemma}[theorem]{Lemma}

\newproclaim{definition}[theorem]{Definition}
\newproclaim{assumption}[theorem]{Assumption}

\newproclaim{remark}[theorem]{Remark}
\newproclaim{example}[theorem]{Example}

\newtheorem{thmm}[theorem]{Theorem}
\newproclaim{defn}[theorem]{Definition}
\newtheorem{prop}[theorem]{Proposition}
\newtheorem{lem}[theorem]{Lemma}
\newtheorem{lemm}{Lemma}[section]
\newproclaim{rque}[theorem]{Remark}
\newproclaim{assum}[theorem]{Assumption}

\newcommand{\Hess}{\operatorname{Hess}}
\newcommand{\Var}{\operatorname{Var}}
\newcommand{\Ent}{\operatorname{Ent}}
\newcommand{\Ric}{\operatorname{Ric}}

\newcommand{\dd}{\mathrm{d}}

\newcommand{\R}{\mathbb{R}}
\newcommand{\N}{\mathbb{N}}
\newcommand{\Z}{\mathbb{Z}}

\newcommand{\cP}{\mathscr{P}}
\newcommand{\cX}{\mathcal{X}}
\newcommand{\cT}{\mathcal{T}}

\newcommand{\cM}{\mathcal{M}}
\newcommand{\cH}{\mathcal{H}}
\newcommand{\cB}{\mathcal{B}}
\newcommand{\tcB}{\widetilde{\mathcal{B}}}
\newcommand{\cW}{\mathcal{W}}
\newcommand{\cL}{\mathcal{L}}

\newcommand{\cG}{\mathcal{G}}
\newcommand{\cK}{\mathcal{K}}
\newcommand{\cE}{\mathcal{E}}
\newcommand{\cA}{\mathcal{A}}

\makeatother

\begin{document}
\begin{frontmatter}

\title{Entropic Ricci curvature bounds for discrete interacting systems\thanksref{T1}}
\runtitle{Ricci bounds for discrete systems}
\thankstext{T1}{Supported by the German Research Foundation
through the Collaborative Research Center 1060 \emph{The Mathematics of
Emergent Effects} and the Hausdorff Center for Mathematics.}

\begin{aug}
\author[A]{\fnms{Max}~\snm{Fathi}\ead[label=e1]{max.fathi@etu.upmc.fr}}
\and
\author[B]{\fnms{Jan}~\snm{Maas}\corref{}\ead[label=e2]{jan.maas@ist.ac.at}}
\runauthor{M. Fathi and J. Maas}
\affiliation{Universit\'{e} Paris 6 and Institute of Science and Technology Austria}
\address[A]{Laboratoire de Probabilit\'es\\
\quad et Mod\`eles Al\'atoires \\
Universit\'{e} Paris 6\\
4, place Jussieu\\
75005 Paris\\
France\\
\printead{e1}}
\address[B]{Institute of Science\\
\quad and Technology Austria (IST Austria)\\
Am Campus 1\\
3400 Klosterneuburg\\
Austria\\
\printead{e2}}
\end{aug}

%
\received{\smonth{1} \syear{2015}}
%
\revised{\smonth{7} \syear{2015}}

%
\begin{abstract}
We develop a new and systematic method for proving entropic Ricci
curvature lower bounds for Markov chains on discrete sets.
Using different methods, such bounds have recently been obtained in
several examples (e.g., \mbox{1-dimensional} birth and death chains, product
chains, Bernoulli--Laplace models, and random transposition models).
However, a general method to obtain discrete Ricci bounds had been
lacking. Our method covers all of the examples above. In addition, we
obtain new Ricci curvature bounds for zero-range processes on the
complete graph. The method is inspired by recent work of Caputo, Dai
Pra and Posta on discrete functional inequalities.
\end{abstract}

%
\end{frontmatter}

\section{Introduction}\label{sec1}

Ricci curvature lower bounds for Riemannian manifolds play a crucial
role in differential geometry, but also in probability theory and
analysis. In particular, such bounds are known to imply important
properties of Brownian motion on the manifold. For example, positive
curvature bounds imply several useful functional inequalities
(logarithmic Sobolev inequality, Poincar\'e inequality), as well as
exponential rates of convergence to equilibrium. Such results were
first obtained using the celebrated Bakry--\'Emery approach developed
in \cite{BE85}.

\subsection{Discrete entropic Ricci curvature}
In \cite{Maas11} and \cite{EM12} a notion of Ricci curvature bounds for
Markov chains was introduced,
based on convexity properties of the relative entropy along geodesics
in the space of probability measures over the state space, for a
well-chosen metric structure. It is the discrete analogue of the
synthetic notion of Ricci curvature bounds in geodesic metric-measure
spaces obtained in the contributions of Lott--Villani \cite{LV09} and
Sturm \cite{St06}.

The main observation behind the Lott--Villani--Sturm definition of
Ricci curvature bounds is that, on a Riemannian manifold $\mathcal{M}$,
the Ricci curvature is bounded from below by some constant $\kappa\in
\R$ if and only if the Boltzmann--Shannon entropy $\mathcal{H}(\rho) =
\int{\rho\log\rho\operatorname{dvol}}$ is $\kappa$-convex along
geodesics in the $L^2$-Wasserstein space of probability measures on
$\mathcal{M}$. However, this definition is not well adapted to discrete
spaces, because the $L^2$-Wasserstein space over a discrete space does
not contain any geodesic.

To circumvent this issue, \cite{Maas11} and \cite{Mi13} introduced a
different metric $\mathcal{W}$. This metric is built from a Markov
kernel in such a way that the heat flow associated to this kernel is
the gradient flow of the entropy with respect to the distance $\mathcal
{W}$. This is an analogue of the continuous situation, where it has
been proven by Jordan, Kinderlehrer and Otto in \cite{JKO} that the
heat flow on ${\R}^d$ is the gradient flow of the entropy with respect
to the Wasserstein metric. A Markov chain is then said to have Ricci
curvature bounded from below by some constant $\kappa$ if the relative
entropy (with respect to the invariant probability measure of the
chain) is $\kappa$-convex along geodesics for the metric $\mathcal{W}$.

It has been shown in \cite{EM12} that a Ricci curvature lower bound has
significant consequences for the associated Markov chain, such as a
(modified) logarithmic Sobolev inequality, a Talagrand transportation
inequality and a Poincar\'e inequality. Therefore, it is desirable to
obtain sharp Ricci curvature bounds in concrete discrete examples. To
this day, there are very few results of this type.

Mielke \cite{Mi13} obtained Ricci curvature bounds for one-dimensional
birth and death processes, and applied these bounds to discretizations
of 1-dimensional Fokker--Planck equations. Erbar and Maas \cite{EM12}
proved a tensorization principle: the Ricci curvature lower bound for a
product system is the minimum of the Ricci curvature bounds of the
components. This property is crucial for applications to
high-dimensional problems. As a special case, the result yields the
sharp Ricci bound for the discrete hypercube $\{0,1\}^n$. The first
high-dimensional results beyond the product-setting have been recently
obtained by Erbar, Maas and Tetali~\cite{EMT}, who obtained Ricci
curvature bounds for the Bernoulli--Laplace model and for the random
transposition model on the symmetric group.

In spite of this progress, a systematic method for obtaining discrete
Ricci curvature bounds has been lacking. In this paper, we present such
a method, which allows us to obtain curvature bounds for several
interacting particle systems on the complete graph. The method allows
us to recover all of the above-mentioned results, and yields new
curvature bounds for zero-range processes on the complete graph.

\subsection{The Bochner-approach}

To explain our method, let us recall that the discrete Ricci curvature
bounds considered in this paper correspond to bounds on the second
derivative of the entropy along $\cW$-geodesics.
In the continuous setting, on a Riemannian manifold $\cM$, the second
derivative of the entropy along a $2$-Wasserstein geodesic $(\rho_t)$
is formally given by
%
\begin{equation}
\label{eq:second-derivative} \frac{\mathrm{ d}^2}{\mathrm{ d}t^2} \cH(\rho_t) = \int
_\cM- \langle\nabla\psi_t, \nabla\Delta
\psi_t \rangle + \frac{1}2 \Delta \bigl(|\nabla
\psi_t|^2 \bigr) \,\dd\rho_t,
\end{equation}
where $\psi: [0,1] \times\cM\to\R$ satisfies the continuity
equation $\partial_t \rho+ \nabla\cdot(\rho\nabla\psi) = 0$ and
the Hamilton--Jacobi equation $\partial_t \psi+ \frac{1}2 |\nabla\psi
|^2 = 0$.
Note that the integrand above can be reformulated using Bochner's
formula, which yields
%
\begin{equation}
\label{eq:Bochner} - \langle\nabla\psi_t, \nabla\Delta\psi_t
\rangle + \tfrac{1}2 \Delta \bigl(|\nabla\psi_t|^2
\bigr) = \tfrac{1}2\bigl|\mathrm{ D^2} \psi\bigr|^2 + \Ric(
\nabla\psi, \nabla\psi).
\end{equation}
Therefore, it follows at least formally that a lower bound $\Ric\geq
\kappa$ implies $\kappa$-convexity of the entropy along geodesics.

In the discrete case, the second derivative of the entropy along
geodesics is given by an expression $\cB(\rho,\psi)$ which somewhat
resembles the right-hand side of \eqref{eq:second-derivative}, but the
dependence on $\rho$ is more complicated. Unfortunately, there does not
seem to be a suitable discrete analogue of Bochner's formula.

In this paper, we present an approach to get around this difficulty. We
obtain a convenient lower bound for the second derivative of the
entropy, which turns out to be very useful; see Theorem~\ref
{thmm:curv-crit} below. This expression has the same features as the
right-hand side of Bochner's formula \eqref{eq:Bochner}: one of the
terms is nonnegative and contains only second-order (discrete)
derivatives, while the other terms contain only first-order
derivatives. The expression is remarkably flexible, since it allows us
to derive sharp Ricci bounds in the rather diverse examples mentioned above.

The work of Caputo, Dai Pra and Posta \cite{CDP09} which inspired this
paper, fits very naturally into this framework. In order to obtain
modified logarithmic Sobolev inequalities, these authors followed the
Bakry--\'Emery strategy, which consists of computing the second
derivative of the entropy along the ``heat equation'' $\partial_t \rho
= \cL\rho$, where $\cL$ is the generator of a reversible Markov chain.
It can be checked that this second derivative can be expressed as $2\cB
(\rho,\log\rho)$. Therefore, the bounds obtained in \cite{CDP09} and
\cite{DP12}
are a
special case of the bounds on $\cB(\rho,\psi)$, which we obtain for
arbitrary $\psi$. We refer to Section~\ref{sec:Bak-Eme} for more details.

\subsection{Organization and notation}

In Section~\ref{sec:pre}, we collect preliminaries on discrete
transport metrics, the associated Riemannian structure, and the notion
of discrete entropic Ricci curvature. Section~\ref{sec:method} contains
the general criterion for Ricci curvature bounds. Various examples are
studied in Section~\ref{sec:examples}. The \hyperref[app]{Appendix}
states a few useful
properties of the logarithmic mean that are used in the proofs.

Throughout the paper, we use the probabilistic notation
\[
\pi[F] = \pi\bigl[F(\eta)\bigr] = \sum_{\eta\in\cX} F(\eta)
\pi(\eta).
\]

\section{Transport metrics and discrete Ricci curvature}
\label{sec:pre}

In this section, we briefly recall some basic facts on the discrete
transportation metric $\cW$, which plays a crucial role in the paper.
This metric has been introduced in \cite{Maas11} in the setting of
finite Markov chains, and (independently) in \cite{Mie11} for
reaction--diffusion systems. Closely related metrics have been
considered in \cite{CHLZ} in the setting of Fokker--Planck equations on
graphs. Our discussion follows \cite{Maas11} with a slightly different
normalisation convention (see also \cite{EM13}).

We work with a discrete, finite space $\mathcal{X}$ and a Markov
generator $\cL$, acting on functions $\psi: \cX\to\R$ by
\[
\cL\psi(\eta) = \sum_{\tilde\eta\in\cX} Q(\eta, \tilde\eta)
\bigl(\psi(\tilde\eta) - \psi(\eta) \bigr).
\]
The transition rates $Q(\eta, \tilde\eta)$ are nonnegative for all
distinct $\eta, \tilde\eta\in\cX$, and we use the convention that
$Q(\eta, \eta) = 0$ for all $\eta\in\cX$.
We assume that $Q$ is \emph{irreducible}, that is, for all $\eta,
\tilde
{\eta} \in\cX$, there exists a sequence $\{\eta_i\}_{i=0}^n
\subseteq
\cX$ such that $\eta_0 = \eta$, $\eta_n = \tilde{\eta}$, and
$Q(\eta_i,
\eta_{i+1}) > 0$ for all $i = 0, \ldots, n-1$. It is a basic result of
Markov chain theory that there is a unique stationary probability
measure $\pi$ on $\mathcal{X}$, which means that
\[
\sum_{\eta\in\mathcal{X}} \pi(\eta) = 1\quad \mbox {and}\quad \pi(\tilde{
\eta}) = \sum_{{\eta} \in\mathcal{X}} Q(\eta, \tilde{\eta})\pi(\eta).
\]
We shall always assume that $\pi$ is \emph{reversible} for $Q$, that
is, the detailed balance equations
\[
Q(\eta, \tilde{\eta})\pi(\eta) = Q(\tilde{\eta}, \eta)\pi (\tilde{\eta})
\]
hold for all $\eta, \tilde{\eta} \in\mathcal{X}$.
We shall refer to $(\cX,Q,\pi)$ as a Markov triple.

\subsection{The nonlocal transport metric $\cW$}

Let
\[
\cP(\mathcal{X}):= \biggl\{ \rho: \mathcal{X} \to\R_+ \Big\vert  \sum
_{\eta\in\mathcal{X}} \rho(\eta)\pi(\eta) = 1 \biggr\}
\]
be the set of \emph{probability densities} on $\mathcal{X}$ (with
respect to the stationary measure $\pi$). We consider the nonlocal
transport metric $\mathcal{W}$ defined for $\rho_0, \rho_1 \in\cP
(\mathcal{X})$ by
\[
\label{eq:def_w} \mathcal{W}(\rho_0, \rho_1)^2
:=\inf_{\rho, \psi} \biggl\{\frac
{1}{2}\int_0^1{
\sum_{\eta, \tilde{\eta} \in
\mathcal{X}} \bigl(\psi_t(\eta) -
\psi_t(\tilde{\eta })\bigr)^2\hat{\rho}_t(
\eta, \tilde{\eta})Q(\eta, \tilde{\eta})\pi (\eta) \,\dd t} \biggr\},
\]
where the infimum runs over all sufficiently smooth curves $\rho:
[0,1] \to\cP(\mathcal{X})$ and $\psi: [0,1] \to\R^{\mathcal{X}}$
satisfying the \emph{continuity equation}
%
\begin{equation}
\label{eq:cont-eq} \frac{\mathrm {d}}{\mathrm{ d}t}\rho_t(\eta) + \sum
_{\tilde{\eta} \in
\mathcal
{X}} \bigl(\psi_t(\tilde{\eta}) -
\psi_t(\eta)\bigr)\hat{\rho }_t(\eta, \tilde{\eta})Q(
\eta, \tilde{\eta}) = 0 \qquad\forall \eta\in\mathcal{X},
\end{equation}
with boundary conditions $\rho|_{t=0} = \rho_0$ and $\rho|_{t=1} =
\rho
_1$. It has been shown in \cite{Maas11} that $\cW$ defines a distance
on the set of probability densities $\cP(\cX)$. Moreover, this distance
is induced by a Riemannian structure on its interior
\[
\cP_*(\cX) := \bigl\{ \rho\in\cP(\cX): \rho(x) > 0 \mbox{ for all $x \in\cX$}
\bigr\}.
\]

This definition is a natural analogue of the Benamou--Brenier
formulation of the Wasserstein metric, with one crucial additional
feature in the discrete setting: the logarithmic mean
%
\begin{equation}\quad
\label{eq:log-mean} \hat{\rho}(\eta, \tilde{\eta}) := \theta\bigl(\rho(\eta), \rho (
\tilde{\eta })\bigr) \qquad\mbox{with } \theta(r,s) = \int_0^1
r^{1-p} s^p \,\dd p = \frac{r - s}{\log r - \log s}
\end{equation}
is used to define, loosely speaking, the value of $\rho$ along the edge
between $\eta$ and $\tilde\eta$. Note that the integral representation
of $\theta$ holds for all $r, s \geq0$, while the alternative
expression requires $r$ and $s$ to be strictly positive. The main
reason for the appearance of this particular mean is the fact that it
allows one to formulate a discrete chain rule $\hat{\rho}\nabla\log
\rho= \nabla\rho$, where
\[
\nabla\psi(\eta,\tilde\eta) := \psi(\tilde\eta) - \psi(\eta)
\]
denotes the discrete gradient.
This identity takes over the role of the usual chain rule $\nabla\log
\rho= \nabla\rho/ \rho$ that conveniently holds in the continuous setting.

In \cite{Maas11,Mie11}, it is shown that the ``heat equation''
$\partial_t \rho= \cL\rho$ is the gradient flow equation for the
relative entropy
\[
\mathcal{H}(\rho) := \sum_{\eta\in\mathcal{X}} \pi(\eta) \rho(\eta)
\log\rho(\eta),
\]
with respect to the Riemannian structure associated with $\mathcal{W}$.
In this sense, the metric $\cW$ is a natural discrete analogue of the
$L^2$-Wasserstein distance. Moreover, at least in some situations \cite
{GM12}, it can be shown that the metric $\cW$ converges to the
$L^2$-Wasserstein in the continuous limit, in the sense of Gromov--Hausdorff.

Following the Lott--Villani--Sturm approach, \cite{Maas11,EM12} gave
the following definition of lower Ricci curvature boundedness.

\begin{defn}[(Ricci curvature lower boundedness)]\label{def:Ricci}
We say that a Markov triple $(\cX,Q,\pi)$ has Ricci curvature bounded
from below by $\kappa\in\R$ if for any constant speed geodesic
$(\rho
_t)_{t \in[0,1]}$ in $(\cP(\mathcal{X}), \mathcal{W})$ we have
\[
\mathcal{H}(\rho_t) \leq(1-t)\mathcal{H}(\rho_0) + t
\mathcal {H}(\rho _1) - \frac{\kappa}{2}t(1-t)\mathcal{W}(
\rho_0, \rho_1)^2.
\]
\end{defn}
We shall use the notation $\Ric(\cX,Q,\pi) \geq\kappa$ [or briefly
$\Ric(Q) \geq\kappa$].
Since $(\cP(\cX),\cW)$ is a geodesic space, this definition is nontrivial.
This notion of curvature has several interesting properties, obtained
in \cite{EM12}:
\begin{itemize}
\item A Markov chain that has $\operatorname{Ric}(Q) \geq\kappa> 0$
satisfies various functional inequalities, such as the modified
logarithmic Sobolev inequality
\[
\label{eq:mLSI} \cH(\rho) \leq\frac{1}{\alpha} \cE(\rho,\log\rho)\qquad \forall\rho\in
\cP(\cX),
\]
with constant $\alpha= 2\kappa$,
as well as a Poincar\'e inequality
\[
\Var_{\pi}(\psi) \leq\frac{1}{\lambda}\cE(\psi,\psi) \qquad\forall\psi:
\mathcal{X} \to\R,
\]
with constant $\lambda= \kappa$. Here, $\Var_{\pi}(\psi) = \pi
[\psi^2]
- \pi[\psi]^2$, and $\cE$ denotes the discrete Dirichlet form given by
\[
\cE(\varphi, \psi) := \frac{1}2 \pi \biggl[ \sum
_{\tilde{\eta}}\bigl(\varphi(\eta) - \varphi(\tilde{\eta })\bigr) \bigl(
\psi (\eta) - \psi(\tilde{\eta})\bigr)Q(\eta, \tilde{\eta}) \biggr].
\]
These functional inequalities imply the exponential convergence estimates
\[
\cH\bigl(e^{tL}\rho\bigr) \leq e^{-\alpha t} \cH(\rho) \quad\mbox{and}\quad
\bigl\|e^{t L} \psi\bigr\|_{L^2(\cX, \pi)} \leq e^{- \lambda t} \| \psi
\|_{L^2(\cX, \pi)}
\]
along the heat equation $\partial_t \rho= \cL\rho$.

\item Ricci curvature bounds tensorize: if $\Ric(\cX_i, Q_i, \pi_i)
\geq\kappa_i$, then the product chain $(\cX_1 \times\cX_2, Q, \pi_1
\otimes\pi_2)$ defined by
\[
Q \bigl((x_1,x_2), (y_1, y_2)
\bigr) = \cases{ %
Q_1(x_1,y_1),&\quad $\mbox{if $x_2 = y_2$},$
\vspace*{2pt}\cr
Q_2(x_2,y_2), & \quad$\mbox{if $x_1 =
y_1$},$
\vspace*{2pt}\cr
0, & \quad$\mbox{otherwise},$}
\]
satisfies $\Ric(\cX_1 \times\cX_2, Q, \pi_1 \otimes\pi_2) \geq
\min
\{ \kappa_1, \kappa_2\}$.
\end{itemize}

\begin{rque}
Besides the entropic Ricci curvature studied in this paper, several
other notions of discrete Ricci curvature have been studied in the
literature, including Ollivier's coarse Ricci curvature \cite{Oll09},
rough Ricci curvature by Bonciocat and Sturm \cite{BS09} and Bakry--\'Emery-type conditions, for example, \cite{LY10}. We refer to \cite
{Oll12} to a survey covering some of these concepts.
\end{rque}

\subsection{The Riemannian structure induced by $\cW$}

It will be useful to describe the Riemannian structure associated to
$\cW$ in more detail.
Let $\cE:= \{ (x,y) \subseteq\cX\times\cX: Q(x,y) > 0 \}$ be the
set of edges in the graph induced by $Q$, and let $\cG$ be the set of
discrete gradients, that is, all functions $\Psi: \cE\to\R$ of the
form $\Psi(\eta, \tilde\eta) = \nabla\psi(\eta, \tilde\eta) :=
\psi
(\tilde\eta) - \psi(\eta)$ for some function $\psi: \cX\to\R$.

Note that, at each $\rho\in\cP_*(\cX)$, the tangent space of $\cP
_*(\cX)$ is
naturally given by $\cT:= \{ \sigma: \cX\to\R| \sum_{\eta\in
\cX} \sigma(\eta)\pi(\eta) = 0\}$.
It can be proved (see \cite{Maas11}, Section~3) that, for each $\rho
\in\cP_*(\cX)$, the mapping
\[
\cK_\rho: \nabla\psi\mapsto\sum_{\tilde{\eta} \in\mathcal
{X}}
\bigl(\psi(\tilde{\eta}) - \psi(\eta)\bigr)\hat{\rho}(\eta , \tilde{\eta})Q(\eta,
\tilde{\eta})
\]
defines an bijection between $\cG$ and $\cT$. At each $\rho\in\cP
_*(\cX)$,
this map allows us to identify the tangent space with $\cG$. In other
words, the continuity equation \eqref{eq:cont-eq} provides an
identification between the ``vertical tangent vector'' $\frac{\mathrm{d}}{\mathrm{ d}t}\rho$ and the ``horizontal tangent vector'' $\nabla\psi$.
A Riemannian structure is then defined by the $\rho$-dependent scalar
product $\langle{\cdot, \cdot}\rangle_\rho: \cG\times\cG\to\R
$ given by
\[
\langle\Phi, \Psi\rangle_{\rho} := \frac{1}{2}\pi \biggl[\sum
_{\tilde
{\eta} \in\mathcal{X}} \Phi(\eta, \tilde{\eta})\Psi(\eta, \tilde{\eta
})\hat{\rho}(\eta, \tilde{\eta})Q(\eta, \tilde{\eta}) \biggr]
\]
for $\Phi, \Psi: \cE\to\R$.
We shall frequently use the notation
\[
\mathcal{A}(\rho, \psi) := \Vert\nabla\psi\Vert_{\rho}^2,
\]
where $\|\cdot\|_\rho$ denotes the norm induced by the scalar product
$\langle{\cdot, \cdot}\rangle_\rho$.

The geodesic equations for this Riemannian structure are given by
\[
\label{eq:geod-equs} %
\cases{ \displaystyle\partial_s
\rho_s(\eta) + \sum_{\tilde\eta\in\cX} \bigl(
\psi_s(\tilde\eta) - \psi_s(\eta) \bigr) \widehat{
\rho_s}(\eta,\tilde\eta) Q(\eta,\tilde\eta) = 0 ,\vspace*{2pt}
\cr
\displaystyle\partial_s \psi_s(\eta) + \frac{1}2 \sum
_{\tilde\eta\in\cX} \bigl( \psi_s(\eta) -
\psi_s(\tilde\eta ) \bigr)^2 (\widehat{
\rho_s})_1(\eta,\tilde{\eta}) Q(\eta,\tilde\eta) = 0,}
\]
where, by abuse of notation, $\hat{\rho}_i$ denotes the partial
derivative of the logarithmic mean defined in \eqref{eq:log-mean}, that
is, for $i =1, 2$,
\[
\label{eq:rho-i} \hat{\rho}_i(\eta,\tilde{\eta}) :=
\partial_i \theta\bigl(\rho(\eta ), \rho (\tilde{\eta})\bigr).
\]
These equations can be regarded as discrete analogues of the continuity
equation and the Hamilton--Jacobi equation, respectively. The Hessian
of the entropy can then be calculated using the expression
\[
\bigl\langle\Hess\mathcal{H}(\rho)\nabla\psi, \nabla\psi\bigr\rangle
_{\rho} = \frac{\mathrm{ d}^2}{\mathrm{ d} s^2} \bigg|_{s = 0} \Ent(\rho_s),
\]
where $(\rho_s)$ solves the geodesic equations with initial conditions
$\rho|_{s= 0} = \rho_0$ and $\psi|_{s= 0} = \psi_0$.
The following result provides an explicit formula for the Hessian which
resembles its continuous counterpart \eqref{eq:second-derivative}; we
refer to \cite{EM12} for a proof. With a slight abuse of notation, we write
\[
\langle\Phi, \Psi\rangle_{\pi} = \frac{1}{2}\pi \biggl[\sum
_{\tilde
{\eta} \in\mathcal{X}} \Phi(\eta, \tilde{\eta})\Psi(\eta, \tilde{\eta
}) Q(\eta, \tilde{\eta}) \biggr].
\]

\begin{prop}[(Identification of the Hessian)] \label{prop:B} For any
$\rho\in\cP_*(\cX)$ and $\psi: \mathcal{X} \to\R$, we have
\begin{eqnarray*}
\cB(\rho, \psi) &:= &\bigl\langle\Hess\mathcal{H}(\rho)\nabla\psi, \nabla \psi
\bigr\rangle_{\rho}
\\
&=& \frac{1}{2}\langle\hat{\cL} \rho\cdot\nabla\psi, \nabla\psi
\rangle_{\pi} - \langle\hat{\rho} \cdot\nabla\psi, \nabla\cL \psi
\rangle_{\pi},
\end{eqnarray*}
where
\[
\hat{\cL}\rho(\eta, \tilde{\eta}) := \widehat{\rho}_1(\eta, \tilde{
\eta })\cL\rho(\eta) + \widehat{\rho}_2(\eta, \tilde{\eta})\cL\rho (
\tilde {\eta}).
\]
\end{prop}

Since the metric structure we consider is Riemannian in $\cP_*(\cX)$, one
expects that $\kappa$-geodesic convexity of $\mathcal{H}$ is
equivalent to
\[
\bigl\langle\Hess\mathcal{H}(\rho)\nabla\psi, \nabla\psi\bigr\rangle
_{\rho} \geq\kappa\Vert\nabla\psi\Vert_{\rho}^2
\]
for any function $\psi$ and any probability density $\rho\in\cP
_*(\cX)$. As
the Riemannian structure is degenerate on the boundary of $\cP
(\mathcal
{X})$, this is not an immediate result, but this issue has been solved
in \cite{EM12}, Theorem~4.4. We thus have the following result.

\begin{prop}[(Characterization of Ricci boundedness)]\label{prop:Ric-char}
Let $\kappa\in\R$. A~Markov triple $(\cX, Q, \pi)$ satisfies
$\operatorname{Ric}(\cX, Q, \pi) \geq\kappa$ if and only if for all
$\rho\in\cP_*(\cX)$ and all functions $\psi: \mathcal{X} \to\R
$, we have
%
\begin{equation}
\label{eq:B-A} \cB(\rho, \psi) \geq\kappa\mathcal{A}(\rho, \psi).
\end{equation}
\end{prop}

In this paper, we provide a systematic method to prove such an
inequality in concrete situations. A very useful criterion will be
given in Section~\ref{sec:method}.

\subsection{The Bakry--\'Emery approach to the logarithmic Sobolev
inequality} \label{sec:Bak-Eme}

Let us now compare our method with the Bakry--\'Emery approach to the
(modified) logarithmic Sobolev inequality (MLSI), which has been
developed in the discrete setting in \cite{BCDP06} and \cite{CDP09}. In
order to prove the MLSI $\cH(\rho) \leq\frac{1}{\alpha} \cE(\rho
,\log
\rho)$, in the Bakry--\'Emery method one considers the behaviour of the
entropy $h(t) := \cH(\rho_t)$ for solutions to the heat equation
$\partial_t \rho_t = \cL\rho_t$. Since $h'(t) = -\cE(\rho_t,\log
\rho
_t)$, the MLSI asserts that
$
h(t) \leq-\frac{1}{\alpha} h'(t)
$.
Rather than approaching this inequality directly, the idea is to
investigate the second derivative of the entropy. Suppose that one
could prove the ``convex entropy decay inequality'' $h''(t) \geq-
{\alpha} h'(t)$. Since $h(t), h'(t) \to0$ as $t \to\infty$, this
inequality implies the MLSI after integration.
Taking into account that
\[
h''(t) = \pi[ \cL\rho_t \cL\log
\rho_t] + \pi \biggl[\frac{(\cL\rho_t)^2}{\rho_t} \biggr],
\]
the following lemma (taken from \cite{CDP09}) summarises this discussion.

\begin{lemma}\label{lem:Bak-Eme}
Let $\alpha> 0$ and suppose that the convex entropy decay inequality
\[
\label{eq:heat-flow-convexity} \pi[ \cL\rho\cL\log\rho] + \pi \biggl[\frac{(\cL\rho)^2}{\rho} \biggr]
\geq {\alpha} \cE(\rho, \log\rho)
\]
holds for all $\rho\in\cP_*(\cX)$.
Then the modified logarithmic Sobolev inequality $\cH(\rho) \leq
\frac
{1}{\alpha} \cE(\rho,\log\rho)$ holds as well. Moreover, \eqref
{eq:heat-flow-convexity} implies the exponential convergence bound
\[
\cE\bigl(e^{t\cL}\rho, \log e^{t\cL}\rho\bigr) \leq
e^{-\alpha t}\cE(\rho, \log\rho)
\]
for all $\rho\in\cP_*(\cX)$.
\end{lemma}
The last assertion follows readily from Gronwall's lemma and is
actually equivalent to \eqref{eq:heat-flow-convexity}. A direct
calculation shows that
\begin{eqnarray*}
\cA(\rho, \log\rho) & =& \cE(\rho, \log\rho),
\\
\cB(\rho, \log\rho) & =& \frac{1}2\pi[ \cL\rho\cL\log\rho] +
\frac{1}2 \pi \biggl[\frac{(\cL\rho)^2}{\rho} \biggr],
\end{eqnarray*}
hence, in view of Proposition~\ref{prop:Ric-char}, the following lemma
follows immediately.

\begin{lemma}\label{lem:Ric-convex decay}
Let $\kappa\in\R$ and suppose that $\Ric(\cX, Q, \pi) \geq\kappa$.
Then the convex entropy decay inequality \eqref{eq:heat-flow-convexity}
holds with $\alpha= 2 \kappa$.
\end{lemma}

The condition $\Ric(\cX, Q, \pi) \geq\kappa$ is in principle strictly
stronger than the convex entropy decay inequality \eqref
{eq:heat-flow-convexity}: we need to check that the inequality $\cB
(\rho
, \psi) \geq\kappa\mathcal{A}(\rho, \psi)$ holds for all functions
$\psi$, and not just for $\psi= \log\rho$. There seems to be no
reason why $\psi= \log\rho$ should always be the extremal case in
this inequality, but we are unaware of an explicit example which
demonstrates this.

In general, the constant we shall obtain here for the lower bound on
Ricci curvature is worse (i.e., smaller) than the best known constant
for the modified logarithmic Sobolev inequality. This was to be
expected, since even in the continuous setting lower bounds on the
Ricci curvature have no reason to yield the optimal constant for the
logarithmic Sobolev inequality in general.

\section{The Bochner method}\label{sec:method}

To explain and apply our method, it will be convenient to use the
following point of view. Given an irreducible and reversible Markov
triple $(\cX,Q,\pi)$, we consider a set $G$ of maps of $\cX$ onto
itself, that represents the set of allowed moves, and a function $c :
\cX\times G \to\R_+$, that represents the jump rates.

\begin{defn}\label{def:mapping}
The pair $(G, c)$ is called a mapping representation of $\cL$ if the
following properties are satisfied:
\begin{longlist}[1.]
\item[1.] The generator $\cL$ can be written as
\[
\cL\psi(\eta) = \sum_{\delta\in G} \nabla _{\delta}
\psi(\eta) c(\eta, \delta),
\]
where $\nabla_{\delta}\psi(\eta) := \psi(\delta\eta) - \psi
(\eta)$.

\item[2.] For any $\delta\in G$, there exists a unique $\delta^{-1} \in G$
satisfying $\delta^{-1}(\delta\eta) = \eta$ for all $\eta$ with
$c(\eta
, \delta) > 0$.

\item[3.] For every $F : \mathcal{X} \times G \to\R$, we have
%
\begin{equation}
\label{eq:reversible}\pi \biggl[ \sum_{\delta\in G} F(\eta, \delta)
c(\eta, \delta) \biggr] = \pi \biggl[ \sum_{\delta\in G} F
\bigl(\delta \eta, \delta^{-1}\bigr) c(\eta, \delta) \biggr].
\end{equation}
\end{longlist}
\end{defn}

Every irreducible reversible Markov chain has a mapping representation.
It is always possible to explicitly build one, by considering the set
of bijections $t_{\eta,\tilde{\eta}} : \mathcal{X} \to\mathcal{X}$
that exchanges $\eta$ and $\tilde{\eta}$, and leaves all other points
unchanged (see \cite{EM12}), but in general it is more natural to work
with a different, smaller set $G$, as we shall see in the examples of
Section~\ref{sec:examples}. Property~3 expresses the reversibility of
the Markov chain.

With this notation, we can rewrite the action functional $\mathcal{A}$ as
\[
\mathcal{A}(\rho, \psi) = \frac{1}{2}\pi \biggl[ \sum
_{\delta\in G} c(\eta, \delta) \bigl(\nabla_{\delta} \psi(\eta)
\bigr)^2\hat{\rho }(\eta, \delta\eta) \biggr].
\]

As discussed in the \hyperref[sec1]{Introduction}, a key element of
the optimal
transport approach to curvature in the continuous setting is the
Bochner identity $\frac{1}{2}\Delta(|\nabla\psi|^2) - \langle
{\nabla \psi , \nabla\Delta\psi}\rangle = |D^2\psi|^2 +
\operatorname{Ric}(\nabla\psi,
\nabla\psi)$. This identity allows one to reformulate
\[
B(\rho, \psi) := \int_{\mathcal{M}}{ \frac{1}{2}\Delta\bigl(|
\nabla \psi|^2\bigr) - \langle\nabla\psi, \nabla\Delta\psi\rangle\,\dd
\rho},
\]
the continuous analogue of $\cB$,
in the equivalent form
\[
B(\rho, \psi) = \int_{\mathcal{M}}{\bigl|D^2
\psi\bigr|^2 + \operatorname {Ric}(\nabla\psi, \nabla\psi)\,\dd\rho}.
\]
Since the first term in this expression is nonnegative, lower bounds
on the Ricci curvature immediately translate into lower bounds on
$B(\rho, \psi)$.

In the discrete setting, as far as we know, there is no direct analogue
of Bochner's identity. To get around this issue (in the context of the
Bakry--\'Emery approach for functional inequalities), it was suggested in
\cite{BCDP06} and \cite{CDP09} to introduce a function $R$ that
satisfies some properties of spatial invariance, in such a way that a
one-sided Bochner inequality holds.

\begin{assum}\label{ass:R}
There exists a function $R : \cX\times G \times G \to\R_+$ such that
\begin{longlist}[(A1)]
\item[(A1)]
$R(\eta, \gamma, \delta) = R(\eta, \delta, \gamma)$ for all $\eta
\in
\cX$ and $\gamma, \delta\in G$;

\item[(A2)]
$ \pi [ \sum_{\gamma, \delta} R(\eta,
\gamma,
\delta)\psi(\eta, \gamma, \delta)  ] = \pi [ \sum_{\gamma, \delta
} R(\eta, \gamma, \delta)\psi(\gamma\eta, \gamma^{-1}, \delta)
 ]$
for all bounded functions $\psi: \mathcal{X} \times G \times G \to\R$;

\item[(A3)]
$\gamma\delta\eta= \delta\gamma\eta$ for all
$\eta\in\mathcal{X}$ and $\gamma, \delta\in G$ with $R(\eta,
\gamma,
\delta) > 0$.
\end{longlist}
\end{assum}

Clearly, (A1) and (A2) imply that for all bounded functions $\psi$,
%
\begin{equation}
\label{eq:R-property} \pi \biggl[ \sum_{\gamma, \delta} R(\eta, \gamma,
\delta)\psi (\eta, \gamma , \delta) \biggr] = \pi \biggl[ \sum
_{\gamma, \delta} R(\eta, \gamma , \delta )\psi\bigl(\delta\eta, \gamma,
\delta^{-1}\bigr) \biggr].
\end{equation}
We will show that a function $R$ that satisfies these assumptions
automatically satisfies the following identity.

\begin{lem} \label{lem:bochner_identity}
Let $\varphi$ and $\psi$ be two real-valued functions on $\mathcal{X}$,
and let $\alpha: \mathcal{X} \times\mathcal{X} \to\R$ be a symmetric
function. Then the following identity holds:
\begin{eqnarray*}
&&\pi \biggl[\sum_{\gamma, \delta}R(\eta, \gamma, \delta)\alpha (
\eta, \delta \eta) \nabla_{\delta}\varphi(\eta) \nabla_{\gamma} \psi(
\eta) \biggr]\\
&&\qquad = \frac{1}{4}\pi \biggl[ \sum_{\gamma, \delta}
R(\eta, \gamma, \delta) \nabla_{\gamma} \bigl[\alpha(\eta, \delta\eta)
\nabla_{\delta} \varphi (\eta) \bigr] \nabla_{\delta}
\nabla_{\gamma} \psi(\eta) \biggr].
\end{eqnarray*}
\end{lem}

In the case $\alpha\equiv1$, the right-hand side in this equation can
be thought of as a discrete analogue of $\int_{\mathcal{M}}{\langle
D^2f, D^2g \rangle \,\dd\rho}$. We can therefore think of this identity as
a sort of weighted discrete Bochner identity.

\begin{rque}
This is a generalization of the identity used in \cite{CDP09}, where
the authors proved this identity in the case $\alpha\equiv1$. In the
applications, we have in mind, we will use this identity with $\alpha$
the logarithmic mean, that is $\alpha(\eta, \delta\eta) = \hat
{\rho
}(\eta, \delta\eta)$.
\end{rque}

\begin{pf*}{Proof of Lemma~\ref{lem:bochner_identity}}
We consider $\eta$, $\delta$ and $\gamma$ such that $R(\eta, \gamma,
\delta) > 0$. Note that, by assumption (iii), $\delta\gamma\eta=
\gamma\delta\eta$. First, we write
\begin{eqnarray*}
&&\nabla_{\gamma} \bigl[\alpha(\eta, \delta\eta)\nabla_{\delta}
\varphi (\eta) \bigr] \nabla_{\delta} \nabla_{\gamma} \psi(\eta) \\
&&\qquad=
\alpha(\gamma\eta, \delta\gamma\eta)\nabla_{\delta} \varphi (\gamma \eta)
\nabla_{\gamma}\psi(\delta\eta) - \alpha(\gamma\eta, \delta\gamma\eta)
\nabla_{\delta} \varphi (\gamma\eta) \nabla_{\gamma}\psi( \eta)
\\
&&\qquad\quad{} - \alpha(\eta, \delta\eta)\nabla_{\delta} \varphi(\eta) \nabla
_{\gamma}\psi(\delta\eta) + \alpha(\eta, \delta\eta)\nabla _{\delta}
\varphi(\eta) \nabla_{\gamma}\psi(\eta).
\end{eqnarray*}

We will show that each of the four terms on the right-hand side of this
equality, when multiplied by $R(\eta, \gamma, \delta)$, summed over all
$\gamma$ and $\delta$ and averaged over $\pi$, yields
\[
\pi \biggl[\sum_{\gamma, \delta} R(\eta, \gamma, \delta)\alpha(
\eta, \delta\eta) \nabla_{\delta}\varphi(\eta) \nabla _{\gamma} \psi(
\eta) \biggr].
\]

For the fourth term, there is nothing to prove. For the third term, we have
\begin{eqnarray*}
&&-\pi \biggl[ \sum_{\gamma, \delta} R(\eta, \gamma, \delta)
\alpha(\eta, \delta\eta)\nabla_{\delta} \varphi(\eta) \nabla _{\gamma}
\psi(\delta\eta) \biggr] \\
&&\qquad= \pi \biggl[ \sum_{\gamma, \delta} R(
\eta, \gamma, \delta) \alpha( \eta, \delta\eta)\nabla_{\delta^{-1}}\varphi (
\delta \eta) \nabla_{\gamma}\psi(\delta\eta) \biggr]
\\
&&\qquad\stackrel{\scriptsize{\eqref{eq:R-property}}} {=} \pi \biggl[ \sum
_{\gamma,
\delta
} R(\eta, \gamma, \delta) \alpha(\delta\eta, \eta )
\nabla_{\delta} \varphi(\eta) \nabla_{\gamma}\psi(\eta) \biggr]\\
&&\qquad = \pi
\biggl[ \sum_{\gamma, \delta} R(\eta, \gamma, \delta) \alpha(
\eta, \delta\eta)\nabla_{\delta} \varphi(\eta) \nabla _{\gamma}\psi(
\eta) \biggr].
\end{eqnarray*}
The same argument works for the second term. Finally, for the first
term, we have
\begin{eqnarray*}
&&\pi \biggl[ \sum_{\gamma, \delta} R(\eta, \gamma ,
\delta)\alpha(\gamma\eta, \delta\gamma\eta)\nabla_{\delta} \varphi (\gamma
\eta) \nabla_{\gamma}\psi(\delta\eta) \biggr]
\\
&&\qquad = -\pi \biggl[ \sum_{\gamma, \delta} R(\eta, \gamma, \delta)
\alpha(\gamma\eta, \delta\gamma\eta)\nabla _{\delta} \varphi(\gamma\eta)
\nabla_{\gamma^{-1}}\psi(\delta\gamma\eta) \biggr]
\\
&&\qquad \stackrel{\mathrm{(A2)}} {=} -\pi \biggl[\sum_{\gamma, \delta} R(\eta,
\gamma, \delta)\alpha(\eta, \delta\eta)\nabla_{\delta} \varphi( \eta)
\nabla_{\gamma}\psi(\delta\eta) \biggr]
\\
&&\qquad= \pi \biggl[ \sum_{\gamma, \delta} R(\eta, \gamma, \delta)
\alpha(\eta, \delta\eta)\nabla_{\delta} \varphi (\eta) \nabla_{\gamma}
\psi(\eta) \biggr], %
\end{eqnarray*}
where the last identity is the one we proved right before for the third
term. This completes the proof.
\end{pf*}

Following \cite{CDP09}, our strategy is to decompose the local weight
$c(\eta, \delta)c(\eta, \gamma)$ as $R(\eta, \delta, \gamma) +
\Gamma
(\eta, \delta, \gamma)$ for a suitable function $\Gamma$. Since the
contribution of $R$ to the curvature is nonnegative, we will only have
to study lower bounds for a functional that depends on $\Gamma$.

\begin{thmm}\label{thmm:curv-crit}
Assume that there exists a function $R$ satisfying Assumption~\ref
{ass:R}, and let
\[
\Gamma(\eta, \gamma, \delta) := c(\eta, \gamma )c(\eta, \delta) - R(\eta,
\gamma, \delta).
\]
Then we have
%
\begin{equation}
\label{eq:B-inequality} \cB(\rho, \psi) \geq (\tcB_1 + \tcB_2 +
\tcB_3) (\rho, \psi),
\end{equation}
where
\begin{eqnarray*}
\tcB_1(\rho, \psi) &:=& \pi \biggl[\sum_{\gamma, \delta}
\Gamma (\eta, \gamma , \delta) \hat{\rho}(\eta, \delta\eta)\nabla_{\delta}
\psi(\eta ) \nabla _{\gamma}\psi(\eta) \biggr],
\\
\tcB_2(\rho, \psi) &:=& \frac{1}{2}\pi \biggl[ \sum
_{\gamma, \delta
} \Gamma (\eta, \gamma, \delta) \bigl(
\nabla_{\delta} \psi(\eta) \bigr)^2 \hat{\rho }_1(
\eta, \delta\eta)\nabla_{\gamma} \rho(\eta) \biggr],
\\
\tcB_3(\rho, \psi) &:=& \frac{1}{4}\pi \biggl[\sum
_{\gamma,
\delta} R(\eta , \gamma, \delta) \hat{\rho}(\eta, \delta\eta)
\bigl(\nabla_{\gamma} \nabla _{\delta} \psi(\eta)\bigr)^2
\biggr]
\end{eqnarray*}
for all $\rho\in\cP(\cX)$ and all $\psi: \cX\to\R$.

As a consequence, if $\tcB_1 + \tcB_2 + \tcB_3 \geq\kappa\mathcal
{A}$, then the Ricci curvature of the Markov chain is bounded from
below by $\kappa$.
\end{thmm}

Inequality \eqref{eq:B-inequality} can be viewed as a discrete
replacement for Bochner's inequality \eqref{eq:Bochner}: the quantity
$\cB(\rho, \psi)$, which is a discrete analogue of the continuous
expression $\int_\cM- \langle\nabla\psi, \nabla\Delta\psi
\rangle
+ \frac{1}2 \Delta (|\nabla\psi|^2 ) \,\dd\rho$, is estimated in
terms of $\tcB_3(\rho, \psi)$, a~nonnegative expression involving only
second-order derivatives, and a number of terms involving first-order
derivatives that may be seen as curvature terms.
%
\begin{rque}
In most of our applications, we will use the trivial bound
$\tcB_3\geq0$,
so that we will only need to obtain a lower bound on the terms $\tcB_1,
\tcB_2$ involving~$\Gamma$.
\end{rque}

\begin{rque}
This result can be used to recover the criterion of Proposition~5.4 in
\cite{EM12}. Indeed, under their assumptions, we can take $R(\eta,
\delta, \gamma) = c(\eta, \delta)c(\eta, \gamma)$ and, therefore,
$\Gamma= 0$. The conclusions easily follow. This criterion has been
applied in \cite{EM12} to obtain the sharp Ricci bound for the discrete
hypercube, as well as nonnegativity of the Ricci curvature for the
discrete cycle $\Z/ N \Z$ for $N \geq2$.
\end{rque}

\begin{pf*}{Proof of Theorem~\ref{thmm:curv-crit}}
The quantity $\mathcal{B}(\rho, \psi)$ from Proposition~\ref{prop:B}
can be written as
\[
\label{eq:decompostion_B_em12} \mathcal{B}(\rho, \psi) = -\langle\hat{\rho} \nabla\psi, \nabla
\cL \psi\rangle_{\pi} + \tfrac{1}{2}\langle\hat{\cL} \rho\nabla
\psi, \nabla\psi\rangle_{\pi},
\]
where
\begin{eqnarray*}
&&\langle\hat{\rho} \nabla\psi, \nabla\cL\psi\rangle_{\pi} \\
&&\qquad=
\frac
{1}{2}\pi \biggl[ \sum_{\gamma, \delta} \nabla
_{\delta} \psi(\eta) \bigl(\nabla_{\gamma} \psi(\delta\eta) c(\delta
\eta, \gamma) - \nabla_{\gamma} \psi(\eta) c(\eta, \gamma) \bigr) c(\eta,
\delta) \hat{\rho}(\eta, \delta\eta) \biggr]
\end{eqnarray*}
and
%
\begin{eqnarray}
\label{eq:lab} \frac{1}{2}\langle\hat{\cL} \rho\nabla\psi,
\nabla\psi\rangle _{\pi
} &= &\frac{1}{4} \pi \biggl[ \sum
_{\gamma, \delta} \bigl(\nabla_{\delta} \psi(\eta)
\bigr)^2 \bigl(\hat{\rho}_1(\eta , \delta \eta)
\nabla_{\gamma}\rho(\eta)c(\eta, \gamma)
\nonumber
\\[-8pt]
\\[-8pt]
\nonumber
&&{} + \hat{\rho}_2(\eta, \delta\eta)\nabla_{\gamma}\rho(\delta
\eta )c(\delta\eta, \gamma) \bigr)c(\eta, \delta) \biggr]. %
\end{eqnarray}
Since the Markov chain is reversible, we can use \eqref{eq:reversible}
to obtain
\begin{eqnarray*}
\label{eq:lab2} &&\pi \biggl[ \sum_{\gamma, \delta}
\nabla_{\delta} \psi(\eta )\nabla _{\gamma} \psi(\delta\eta) c(\delta
\eta, \gamma) c(\eta, \delta ) \hat {\rho}(\eta, \delta\eta) \biggr]
\nonumber
\\[-8pt]
\\[-8pt]
\nonumber
&&\qquad= - \pi
\biggl[ \sum_{\gamma, \delta} \nabla_{\delta} \psi(\eta)
\nabla_{\gamma} \psi(\eta) c(\eta, \gamma) c(\eta, \delta)\hat{\rho}(\eta,
\delta\eta) \biggr], %
\end{eqnarray*}
so that
\begin{eqnarray*}
-\langle\hat{\rho}  \nabla\psi, \nabla\cL\psi\rangle_{\pi} &= &\pi
\biggl[ \sum_{\gamma, \delta} \nabla_{\delta} \psi(\eta)
\nabla _{\gamma} \psi(\eta) c(\eta, \gamma) c(\eta, \delta)\hat{\rho}(\eta,
\delta\eta) \biggr]
\\
&= &\pi \biggl[ \sum_{\gamma, \delta}\Gamma(\eta, \gamma,
\delta) \nabla _{\delta} \psi(\eta)\nabla_{\gamma} \psi(\eta) \hat{
\rho}(\eta , \delta \eta) \biggr]\\
&&{} + \pi \biggl[ \sum_{\gamma, \delta}
R(\eta, \gamma, \delta) \nabla_{\delta} \psi(\eta)\nabla_{\gamma}
\psi(\eta) \hat{\rho }(\eta, \delta\eta) \biggr].
\end{eqnarray*}
The first term on the right-hand side of this equality we leave
unchanged, and for the second we use Lemma~\ref{lem:bochner_identity}
with $\alpha= \hat{\rho}$ to get
%
\begin{eqnarray}
\label{eq:equality1} -\langle\hat{\rho} \nabla\psi, \nabla\cL\psi
\rangle_{\pi} &= &\pi \biggl[ \sum_{\gamma, \delta}
\Gamma(\eta, \gamma, \delta) \nabla _{\delta} \psi(\eta)
\nabla_{\gamma} \psi(\eta) \hat{\rho}(\eta , \delta \eta) \biggr]
\nonumber
\\[-8pt]
\\[-8pt]
\nonumber
&&{} + \frac{1}{4}\pi \biggl[ \sum_{\gamma, \delta} R(\eta,
\gamma, \delta ) \nabla_{\delta}\nabla_{\gamma} \psi(\eta)
\nabla_{\gamma
}\bigl[\hat{\rho }(\eta, \delta\eta) \nabla_{\delta}
\psi(\eta)\bigr] \biggr].
\end{eqnarray}

We now take a look at the second term in \eqref{eq:lab}. Using the
reversibility of the Markov chain and the fact that $\hat{\rho
}_2(\delta\eta, \eta) = \hat{\rho}_1(\eta, \delta\eta)$, we obtain
\begin{eqnarray*}
&&\pi \biggl[ \sum_{\gamma, \delta} \bigl(\nabla_{\delta}
\psi(\eta ) \bigr)^2\hat{\rho}_2(\eta, \delta\eta)
\nabla_{\gamma}\rho(\delta \eta )c(\delta\eta, \gamma) c(\eta, \delta)
\biggr]
\\
&&\qquad = \pi \biggl[ \sum_{\gamma, \delta} \bigl(
\nabla_{\delta} \psi (\eta) \bigr)^2\hat{\rho}_1(
\eta, \delta\eta)\nabla_{\gamma}\rho(\eta )c(\eta, \gamma) c(\eta, \delta)
\biggr],
\end{eqnarray*}
and, therefore,
\begin{eqnarray*}
\frac{1}2\langle\hat{\cL}\rho\nabla\psi, \nabla\psi\rangle_{\pi
}
&= &\frac{1}{2}\pi \biggl[ \sum_{\gamma, \delta} \bigl(
\nabla_{\delta} \psi(\eta ) \bigr)^2\hat{\rho}_1(
\eta, \delta\eta)\nabla_{\gamma}\rho(\eta )c(\eta, \gamma) c(\eta, \delta)
\biggr]
\\
&=& \frac{1}{2}\pi \biggl[ \sum_{\gamma, \delta} \Gamma(
\eta, \gamma, \delta) \bigl(\nabla_{\delta} \psi(\eta) \bigr)^2
\hat{\rho }_1(\eta, \delta \eta)\nabla_{\gamma}\rho(\eta)
\biggr]
\\
&&{} + \frac{1}{2}\pi \biggl[ \sum_{\gamma, \delta} R(\eta,
\gamma, \delta ) \bigl(\nabla_{\delta} \psi(\eta) \bigr)^2\hat{
\rho}_1(\eta, \delta\eta )\nabla_{\gamma}\rho(\eta) \biggr]\\
& =:&
T_1 + T_2.
\end{eqnarray*}
Again, we leave the first term on the right-hand side of this last
equation unchanged. For the second term, we use assumption (A2) to obtain
\[
T_2 = \frac{1}2\pi \biggl[ \sum
_{\gamma, \delta} R(\eta, \gamma, \delta ) \bigl(\nabla_{\delta}
\psi(\eta) \bigr)^2\hat{\rho}_2(\eta, \delta \eta)\nabla
_{\gamma}\rho(\delta\eta) \biggr],
\]
and thus
\[
T_2 = \frac{1}{4}\pi \biggl[ \sum
_{\gamma, \delta} R(\eta, \gamma, \delta ) \bigl(\nabla_{\delta}
\psi(\eta) \bigr)^2 \bigl(\hat{\rho}_1(\eta , \delta \eta)
\nabla_{\gamma}\rho(\eta) + \hat{\rho}_2(\eta, \delta\eta )
\nabla _{\gamma}\rho(\delta\eta) \bigr) \biggr].
\]
From (i) and (ii) of Lemma~\ref{lem:theta} we deduce the inequality
\[
\hat{\rho}_1(\eta, \delta\eta)\nabla_{\gamma}\rho(\eta) + \hat
{\rho }_2(\eta, \delta\eta)\nabla_{\gamma}\rho(\delta\eta) \geq
\nabla_{\gamma
}\hat{\rho}(\eta, \delta\eta).
\]
Therefore, we have
%
\begin{eqnarray}
\label{eq:last} &&\frac{1}{2}\langle\hat{\cL} \rho\nabla\psi, \nabla\psi
\rangle _{\pi
}\nonumber\\
&&\qquad\geq\frac{1}{2}\pi \biggl[ \sum
_{\gamma, \delta} \Gamma(\eta, \gamma, \delta) \bigl(\nabla_{\delta}
\psi(\eta) \bigr)^2\hat{\rho }_1(\eta, \delta \eta)
\nabla_{\gamma}\rho(\eta) \biggr]
\\
&&\qquad\quad{} +\frac{1}{4}\pi \biggl[ \sum_{\gamma, \delta} R(\eta,
\gamma, \delta) \bigl(\nabla_{\delta} \psi(\eta) \bigr)^2
\nabla_{\gamma}\hat{\rho }(\eta, \delta\eta) \biggr].\nonumber
\end{eqnarray}
To deduce our result from \eqref{eq:equality1} and \eqref{eq:last}, all
that is left is to show that
%
\begin{eqnarray}
\label{eq:last_identity_R} S &:=& \pi \biggl[ \sum_{\gamma, \delta}
R(\eta, \gamma, \delta ) \bigl( \bigl(\nabla_{\delta} \psi(\eta)
\bigr)^2\nabla_{\gamma}\hat{\rho }(\eta, \delta\eta)\nonumber \\
&&{}+
\nabla_{\delta}\nabla_{\gamma} \psi(\eta) \nabla_{\gamma} \bigl[
\hat {\rho}(\eta, \delta\eta) \nabla_{\delta} \psi(\eta) \bigr] \bigr)
\biggr]
\\
& =&\pi \biggl[ \sum_{\gamma, \delta} R(\eta, \gamma, \delta)
\hat{\rho }(\eta, \delta\eta) \bigl(\nabla_{\gamma} \nabla_{\delta}
\psi(\eta )\bigr)^2 \biggr].\nonumber
\end{eqnarray}
We observe that
\begin{eqnarray*}
&& \bigl(\nabla_{\delta}\psi(\eta) \bigr)^2
\nabla_{\gamma}\hat{\rho }(\eta, \delta\eta) + \nabla_{\delta}
\nabla_{\gamma} \psi(\eta) \nabla_{\gamma} \bigl[\hat{\rho}(\eta,
\delta\eta)\nabla_{\delta
}\psi (\eta) \bigr]
\\
& &\qquad= \hat\rho(\gamma\eta, \delta\gamma\eta) \bigl( \nabla_\gamma
\nabla_\delta\psi(\eta) \bigr)^2\\
&&\qquad\quad{} - \hat\rho(\eta, \delta\eta)
\nabla_\delta\psi(\gamma\eta) \nabla_\delta\psi(\eta) + \hat
\rho(\gamma\eta, \delta\gamma\eta) \nabla_\delta\psi(\gamma\eta)
\nabla_\delta\psi(\eta) ,
\end{eqnarray*}
for all $\eta\in\cX$ and $\gamma, \delta\in G$ with $\delta\gamma
\eta= \gamma\delta\eta$. In view of (A3), we may insert this identity
in the expression for $S$, which yields
\begin{eqnarray*}
S &=&\pi \biggl[ \sum_{\gamma, \delta} R(\eta, \gamma, \delta)
\hat {\rho }(\gamma\eta, \delta\gamma\eta) \bigl(\nabla_{\gamma}
\nabla_{\delta} \psi (\eta)\bigr)^2 \biggr]\\
&&{} - \pi \biggl[ \sum
_{\gamma, \delta} R(\eta, \gamma, \delta) \hat {\rho }(\eta,
\delta\eta)\nabla_{\delta}\psi(\gamma\eta) \nabla _{\delta}\psi (
\eta) \biggr]
\\
&&{} + \pi \biggl[ \sum_{\gamma, \delta} R(\eta, \gamma, \delta)
\hat{\rho }(\gamma\eta, \delta\gamma\eta)\nabla_{\delta}\psi(\gamma\eta )
\nabla _{\delta}\psi( \eta) \biggr].
\end{eqnarray*}
It follows from (A2) that the second and third term in the latter
formula cancel each other. Another application of (A2) shows that the
first term equals the right-hand side of \eqref{eq:last_identity_R},
which was the last element we needed to complete the proof.
\end{pf*}

\section{Examples}
\label{sec:examples}

\subsection{Birth and death processes}

In this section, we consider the case of birth and death processes on
$\N= \{0,1, \ldots\}$. These are Markov chains with generator
\[
\mathcal{L}\psi(n) = a(n) \bigl(\psi(n+1) - \psi(n) \bigr) + b(n) \bigl(\psi
(n-1) - \psi(n) \bigr),
\]
where $a$ and $b$ are nonnegative functions on $\N$, such that $b(0) =
0$. The set of allowed moves is $G := \{+, -\}$, where $+(n) = n+1$ and
$-(n) = (n-1)\mathbh{1}_{\{n>0\}}$. Following the notation of the
previous section, we write $\nabla_{\pm}\psi(n) = \psi(n \pm1)
-\psi
(n)$. The generator can therefore be written as
\[
\mathcal{L}\psi(n) = a(n)\nabla_+ \psi(n) + b(n)\nabla_- \psi(n),
\]
and, in accordance with our notation, we set $c(n,+) = a(n)$ and
$c(n,-) = b(n)$.

We assume that this Markov chain is irreducible, and that there exists
a probability measure $\pi$ on $\N$ satisfying the detailed balance condition
\[
a(n)\pi(n) = b(n+1)\pi(n+1).
\]
Clearly, the latter condition is satisfied if and only if
\[
\sum_{n=1}^\infty\frac{a(n-1) \cdots a(0)}{b(n)
\cdots
b(1)} < \infty.
\]
When applying our method to study curvature bounds for this Markov
chain, we obtain the following result.

\begin{thmm}\label{th4.1}
Let $\kappa\in\R_+$. Assume that the rate of birth $a$ is
nonincreasing, and that the rate of death $b$ is nondecreasing.
Assume moreover that
%
\begin{eqnarray}
\label{eq:assumpt_rates_bd} &&\tfrac{1}2 \bigl(a(n) - a(n+1) + b(n+1)
- b(n) \bigr)
\nonumber
\\[-8pt]
\\[-8pt]
\nonumber
&&\qquad{}+
\tfrac{1}{2}\Theta \bigl(a(n) - a(n+1), b(n+1) - b(n) \bigr) \geq\kappa
\end{eqnarray}
for all $n \in\N$, where
\[
\Theta(\alpha, \beta) = \inf_{s, t > 0} \theta (s,t) \biggl(
\frac{\alpha}{s} + \frac{\beta}{t} \biggr).
\]
Under these assumptions, the birth and death process has Ricci
curvature bounded from below by $\kappa$.
\end{thmm}

\begin{rque}\label{rem:countable}
Strictly speaking, to include this example in our framework we need to
assume that the transition rates $a(n), b(n)$ vanish whenever $n$ is
sufficiently large, so that the state spaces becomes finite. However,
it is reasonable to expect that the result below and its proof remain
valid without this assumption. A rigorous inclusion of the countable
setting would require us to modify some technical arguments from \cite
{Maas11,EM12}, which is beyond the scope of the present paper.
\end{rque}

\begin{rque}
\label{rem:crit}
The same criterion has been obtained by Mielke \cite{Mi13} using a
different method based on diagonal dominance of the matrix representing
the Hessian of the entropy. Mielke worked on a finite state space $\{0,
1, \ldots, N\}$ rather than on $\N$, but this does not effect the argument.
\end{rque}

\begin{pf*}{Proof of Theorem \ref{th4.1}}
First of all, using the reversibility condition, we can rewrite the
action functional as
\begin{eqnarray*}
\mathcal{A}(\rho, \psi) &= &\frac{1}{2}\pi \bigl[a(n) \bigl(\nabla_+\psi
(n)\bigr)^2\hat{\rho}(n, n+1) + b(n) \bigl(\nabla_-\psi(n)
\bigr)^2\hat{\rho}(n, n-1) \bigr]
\\
&=& \pi \bigl[a(n) \bigl(\nabla_+\psi(n)\bigr)^2\hat{\rho}(n, n+1)
\bigr].
\end{eqnarray*}
Following \cite{CDP09}, we define the function $R : \N\times{\{+,-\}
}^2 \to\R$ by
\[
\cases{ R(n, +, +) = a(n)a(n+1), \vspace*{2pt}
\cr
R(n, -, -) =
b(n)b(n-1), \vspace*{2pt}
\cr
R(n, +, -) = R(n, -, +) = a(n)b(n)} %
\]
so that $\Gamma(n, +, +) = - a(n)\nabla_+ a(n)$, $\Gamma(n, -, -) = -
b(n)\nabla_- b(n)$, and $\Gamma(n, +,   -) = \Gamma(n, -, +) = 0$.
It is not hard to check that Assumption~\ref{ass:R} is satisfied. In
this case, we have
%
\begin{eqnarray}
\label{eq:equality3} \tcB_1(\rho,\psi)& =& \pi \biggl[ \sum
_{\gamma, \delta} \Gamma (\eta, \gamma , \delta) \hat{\rho}(\eta,
\delta\eta)\nabla_{\delta}\psi(\eta ) \nabla _{\gamma}\psi(\eta)
\biggr]
\nonumber
\\
&=& - \pi \bigl[a(n)\nabla_+ a(n) \bigl(\nabla_+\psi(n) \bigr)^2\hat
{\rho}(n, n+1)
\nonumber
\\[-8pt]
\\[-8pt]
\nonumber
&&{}+ b(n)\nabla_-b(n) \bigl(\nabla_-\psi(n)\bigr)^2\hat{
\rho}(n, n-1) \bigr]
\\
&=& \pi \bigl[a(n) \bigl( \nabla_+ b(n) - \nabla_+ a(n) \bigr) \bigl(\nabla _+
\psi(n) \bigr)^2\hat{\rho}(n, n+1) \bigr].
\nonumber
\end{eqnarray}
Moreover, using reversibility,
\begin{eqnarray*}
\tcB_2(\rho,\psi) & =& \frac{1}{2}\pi \biggl[ \sum
_{\gamma, \delta
} \Gamma (\eta, \gamma, \delta) \bigl(
\nabla_{\delta} \psi(\eta) \bigr)^2 \hat{\rho }_1(
\eta, \delta\eta)\nabla_{\gamma} \rho(\eta) \biggr]
\\
& = &- \frac{1}{2}\pi \bigl[a(n) \nabla_+a(n) \bigl(\nabla_+\psi (n)
\bigr)^2\hat{\rho}_1(n, n+1)\nabla_+\rho(n)
\\
&&{} + b(n)\nabla_- b(n) \bigl(\nabla_-\psi(n) \bigr)^2\hat{
\rho}_1(n, n-1)\nabla_-\rho(n) \bigr]
\\
&= &-\frac{1}{2}\pi \bigl[a(n)\nabla_+ a(n) \bigl(\nabla_+\psi (n)
\bigr)^2\hat{\rho}_1(n, n+1)\nabla_+\rho(n)
\\
& &{}+ a(n)\nabla_+ b(n) \bigl(\nabla_+\psi(n) \bigr)^2\hat{
\rho}_2(n, n+1)\nabla_+\rho(n)\bigr ].
\end{eqnarray*}
A direct computation shows that the partial derivatives $\theta_i(s,t)
= \partial_i \theta(s,t)$ satisfy the identities
\[
(t - s)\theta_1(s, t) = \frac{\theta(s,t)^2}{s} - \theta(s,t),\qquad  (s - t)
\theta_2(s, t) = \frac{\theta(s,t)^2}{t} - \theta(s,t),
\]
for any $s, t > 0$.
Therefore, recalling that $\hat{\rho}(n, n+1) = \theta (\rho
(n), \rho
(n+1) )$, we have
\begin{eqnarray}
\label{eq:tilde-B} \tcB_2(\rho,\psi) &=&-\frac{1}{2}\pi
\bigl[a(n) \bigl(\nabla_+ b(n) - \nabla_+ a(n) \bigr) \bigl(\nabla_+\psi(n)
\bigr)^2\hat{\rho}(n, n+1) \bigr]
\nonumber
\\[-8pt]
\\[-8pt]
\nonumber
&&{} + \frac{1}{2}\pi \biggl[a(n) \biggl( \frac{\nabla_+ b(n)}{\rho
(n+1)}-
\frac{\nabla_+a(n)}{\rho(n)} \biggr) \bigl(\nabla_+\psi (n) \bigr)^2\hat{
\rho}(n, n+1)^2 \biggr].
\nonumber
\end{eqnarray}
Using the definition of $\Theta$, we obtain the inequality
\[
\frac{\nabla_+b(n)}{\rho(n+1)} -\frac{\nabla_+a(n)}{\rho(n)} \geq\frac{\Theta ({-}\nabla_+ a(n), \nabla_+b(n) )}{\hat
{\rho
}(n,n+1)}.
\]
Summing \eqref{eq:equality3} and \eqref{eq:tilde-B}, and substituting
the latter inequality, it follows from the assumption \eqref
{eq:assumpt_rates_bd} that
\[
(\tcB_1+\tcB_2) (\rho, \psi) \geq\kappa\mathcal{A}(\rho,
\psi).
\]
The result follows by an application of Theorem~\ref{thmm:curv-crit}.
\end{pf*}

\subsection{Zero-range processes}

We now consider interacting particle systems on the complete graph,
whose sites are labeled $1, \ldots,L$, where $L$ is the number of
sites. Configurations of particles are given by collections of
nonnegative integers $\eta_1,\ldots, \eta_L$, where $\eta_x$ denotes
the number of particles at site $x$. The whole configuration will be
denoted by $\eta\in S := \N^L$.

Changes in the configuration are made by moving a particle (if any)
from a site $x$ to a different site $y$. Allowed moves are therefore
given by maps of the form $\eta\mapsto\eta^{xy}$, where $x \neq y$
are two different sites, and $\eta^{xy}$ is given by
\[
\cases{ \eta^{xy}_z = \eta_z, &\quad $
\mbox{if } z \notin\{x,y\} \mbox{ or } \eta_x = 0 $, \vspace*{2pt}
\cr
\eta^{xy}_x = \eta_x - 1, &\quad $\mbox{if }
\eta_x > 0 $, \vspace *{2pt}
\cr
\eta^{xy}_y =
\eta_y + 1, &\quad $\mbox{if } \eta_x > 0 $.} %
\]
The set of allowed moves is $G := \{xy; x \neq y \}$,
where $xy$ denotes the map $\eta\mapsto\eta^{xy}$. We write $\nabla
_{xy}$ to denote the corresponding discrete gradient, and note that
$(xy)^{-1} = (yx)$ in the sense of Definition~\ref{def:mapping}.

In this section, we consider zero-range processes on the complete graph
with $L$ vertices. In such processes, the jump rate from site $x$ to
site $y$ depends only on $x$ and on the number of particles present at
$x$. The rates are thus given by a family of functions $c_x : \N\to
[0, \infty)$ such that $c_x(0) = 0$ and $c_x(n) > 0$ for all $n > 0$
and $x \in\{1, \ldots, L\}$. Here, $c_x(\eta_x)$ is the rate at which
a particle is moved from site $x$ to a site $y$, with $y$ chosen
randomly, with uniform probability on $\{1, \ldots,L\}$.
The generator of the Markov chain we just described is
\[
\mathcal{L}f(\eta) = \frac{1}{L}\sum_{x,y}
c_x(\eta_x)\nabla_{xy}f(\eta).
\]
This dynamic conserves the total number of particles $N := \sum_x \eta
_x$. We define a probability measure $\pi_N$ on configurations with $N$
particles as
\[
\pi_N(\eta) := \frac{1}{Z_N}\prod
_{x=1}^{{L}} \prod_{k=1}^{\eta
_x}
\frac
{1}{c_x(k)},
\]
with the usual convention that the latter product equals $1$ whenever
$\eta_x = 0$.
The constant $Z_N$ is a normalization constant, which is indeed finite,
since there are a finite number of configurations with $N$ particles.
It can then be checked that the Markov chain is reversible with respect
to $\pi_N$. In the sequel, we shall systematically omit the subscript
$N$, which we consider fixed.

When applying our method to study the curvature of this process, we
obtain the following result.

\begin{thmm}\label{thmm:zrp}
Assume that there exists $c > 0$ and $\delta\in[0, 2c]$ such that
%
\begin{equation}
\label{eq:zrp-ass} c \leq c_x(n+1) - c_x(n) \leq c + \delta\qquad
\mbox{for all $x \in\{1,\ldots, L\}$ and $n \geq0$}.
\end{equation}
Then the Ricci curvature is bounded from below by $\frac{c}{2} - \frac
{5\delta}{4}$. In particular, if $\delta< \frac{2}{5}c$, then we
obtain a positive lower bound on the Ricci curvature.
\end{thmm}

Note that the assumption implies that the rates $c_x$ are strictly
increasing. Of particular interest is the fact that the bound does not
depend on either $L$ or $N$.

\begin{rque}
Simply assuming that the rates are strictly increasing is not enough to
ensure Ricci curvature is positive. In Section~4.2 of \cite{CDP09}, a
simple example of a zero-range process for which $\sup c_x(1) - \inf
c_x(1)$ is large was shown to exhibit nonconvex decay of the entropy
and, therefore, the Ricci curvature of the Markov kernel cannot be
positive. So some assumption of the form $\delta< Kc$ is necessary. On
the other hand, we do not know if the assumption $\delta< 2c/5$ is
optimal. Assuming a uniform positive lower bound on the increase of the
rates is known to be necessary to obtain a lower bound on the curvature
that is independent of the system size and number of particles. Indeed,
it was shown in \cite{Mor} that for the case of constant rates, the
spectral gap is of order $\frac{L^2}{L^2 + N^2}$.
\end{rque}

\begin{rque}
For independent random walks on the complete graph (which corresponds
to $\delta= 0$ and $c = 1$), we recover the size-independent lower
bound on Ricci curvature that was obtained in Example~5.1 of \cite{EM12}.

The estimate on the modified logarithmic Sobolev inequality that can be
deduced from the curvature bound has a worse dependence on $\delta$
than the one obtained in \cite{CDP09} [which would correspond to a
constant $(c-\delta)/2$]. The modified logarithmic Sobolev inequality
for this model has also been studied in \cite{CP07}.
\end{rque}

\begin{pf*}{Proof of Theorem \ref{thmm:zrp}}
First, we observe that
\[
\mathcal{A}(\rho, \psi) = \frac{1}{2L}\pi \biggl[\sum
_{x,y} c_x(\eta_x) \bigl(
\nabla_{xy}\psi(\eta)\bigr)^2\hat{\rho}\bigl(\eta, \eta
^{xy}\bigr) \biggr].
\]
Following \cite{CDP09}, we define the function $R$ as
\[
R(\eta, xy, uv) = \cases{ \displaystyle
\frac{1}{L^2}c_x(
\eta_x)c_u(\eta_u), &\quad $\mbox{for } x \neq u,$
\vspace*{2pt}\cr
\displaystyle\frac{1}{L^2}c_x(\eta_x)c_x(
\eta_x - 1), &\quad $\mbox{otherwise},$ }
\]
with the understanding that $c_x(-1) := 0$,
so that $\Gamma(\eta, xy, uv) = 0$ if $x \neq u$, and $\Gamma(\eta, xy,
xv) = \frac{1}{L^2} c_x(\eta_x) (c_x(\eta_x) - c_x(\eta_x -
1) )$.
It follows that\vspace*{2pt}
%
\begin{eqnarray}
\label{eq:eq5} \tcB_2(\rho, \psi) &=& \frac{1}{2}
\pi \biggl[ \sum_{\gamma, \delta
} \Gamma(\eta, \gamma, \delta)
\bigl(\nabla_{\delta} \psi(\eta ) \bigr)^2 \hat{
\rho}_1(\eta, \delta\eta)\nabla_{\gamma} \rho(\eta) \biggr]
\nonumber
\\[1pt]
&=& \frac{1}{2L^2}\pi \biggl[ \sum_{x,y,v}
c_x(\eta _x) \bigl(c_x(
\eta_x) - c_x(\eta_x - 1)\bigr) \bigl(
\nabla_{xy}\psi(\eta)\bigr)^2
\\[1pt]
&&{}\times \hat{\rho}_1\bigl(
\eta, \eta ^{xy}\bigr)\nabla_{xv}\rho(\eta) \biggr].\nonumber
\end{eqnarray}
Using the reversibility assumption, this quantity is also equal to\vspace*{2pt}
\begin{eqnarray*}
&&\frac{1}{2L^2}\pi \biggl[ \sum_{x,y,v}
c_x(\eta_x) \bigl(c_y(\eta_y
+ 1) - c_y(\eta_y) \bigr) \bigl(\nabla_{xy}
\psi(\eta) \bigr)^2\hat{\rho }_2\bigl(\eta, \eta
^{xy}\bigr)
\\[1pt]
&&\qquad{}\times  \bigl(\rho\bigl(\eta^{xv}\bigr) - \rho\bigl(
\eta^{xy}\bigr) \bigr) \biggr].
\end{eqnarray*}
Adding this identity to \eqref{eq:eq5}, we obtain\vspace*{2pt}
\begin{eqnarray*}
\hspace*{-5pt}&&\tcB_2(\rho, \psi) \\[1pt]
\hspace*{-5pt}&&\qquad= \frac{1}{4L^2}\pi \biggl[ \sum
_{x,y,v} c_x(\eta_x) \bigl(
\nabla _{xy}\psi (\eta) \bigr)^2 \\[1pt]
\hspace*{-5pt}&&\qquad\quad{}\times \bigl(
\bigl(c_x(\eta_x) - c_x(\eta_x
- 1) \bigr)\hat {\rho }_1\bigl(\eta, \eta^{xy}\bigr)
\nabla_{xv}\rho(\eta)\\[1pt]
\hspace*{-5pt}&&\qquad\quad{}+ \bigl(c_y(\eta_y + 1) - c_y(
\eta_y) \bigr)\hat{\rho}_2\bigl(\eta, \eta
^{xy}\bigr) \bigl(\rho\bigl(\eta^{xv}\bigr) - \rho\bigl(
\eta^{xy}\bigr) \bigr) \bigr) \biggr]
\nonumber
\\[1pt]
\hspace*{-5pt}&&\qquad= \frac{1}{4L^2}\pi \biggl[ \sum_{x,y,v}
c_x(\eta_x) \bigl(\nabla _{xy}\psi (\eta)
\bigr)^2\rho\bigl(\eta^{xv}\bigr) \\[1pt]
\hspace*{-5pt}&&\qquad\quad{}\times\bigl(
\bigl(c_x(\eta_x) - c_x(\eta_x
- 1)\bigr)\hat {\rho}_1\bigl(\eta, \eta^{xy}\bigr)
 + \bigl(c_y(\eta_y + 1) - c_y(
\eta_y) \bigr)\hat{\rho}_2\bigl(\eta, \eta
^{xy}\bigr) \bigr) \biggr]
\nonumber
\\[1pt]
\hspace*{-5pt}& &\qquad\quad{}- \frac{1}{4L^2}\pi \biggl[ \sum_{x,y,v}
c_x(\eta_x) \bigl(\nabla _{xy}\psi (\eta)
\bigr)^2 \bigl( \bigl(c_x(\eta_x) -
c_x(\eta_x - 1) \bigr)\hat {\rho }_1\bigl(
\eta, \eta^{xy}\bigr)\rho(\eta)
\\[1pt]
\hspace*{-5pt}&&\qquad\quad{} + \bigl(c_y(\eta_y + 1) - c_y(
\eta_y) \bigr)\hat{\rho}_2\bigl(\eta, \eta
^{xy}\bigr)\rho\bigl(\eta^{xy}\bigr) \bigr) \biggr]
\\[1pt]
\hspace*{-5pt}&&\qquad =: I_1 + I_2. %
\end{eqnarray*}
Taking into account that $\hat{\rho}_i \geq0$, we use \eqref
{eq:zrp-ass} and the assumption that $c \geq\delta/2$ to obtain the
lower bound
\begin{eqnarray*}
I_1 &\geq&\frac{c}{4L^2}\pi \biggl[ \sum
_{x,y,v} c_x(\eta_x) \bigl(
\nabla_{xy}\psi(\eta) \bigr)^2\rho\bigl(\eta^{xv}
\bigr) \bigl(\hat {\rho}_1\bigl(\eta, \eta^{xy}\bigr) +
\hat{\rho}_2\bigl(\eta, \eta^{xy}\bigr) \bigr) \biggr]
\\
&\stackrel{y \leftrightarrow v} {\geq} &\frac{\delta
}{8L^2}\pi \biggl[ \sum
_{x,y,v} c_x(\eta_x) \bigl(
\nabla_{xv}\psi(\eta) \bigr)^2\rho\bigl(\eta^{xy}
\bigr) \bigl(\hat{\rho }_1\bigl(\eta, \eta ^{xv}\bigr) +
\hat{\rho}_2\bigl(\eta, \eta^{xv}\bigr) \bigr) \biggr]
\\
&=:& I_3. %
\end{eqnarray*}
This quantity is nonnegative, but we will need it later to compensate
a term that appears when computing $\tcB_1(\rho, \psi)$.
Using \eqref{eq:zrp-ass} once more and Lemma~\ref{lem:theta}(ii),
we obtain
\begin{eqnarray*}
I_2 &\geq&- \frac{c + \delta}{4L^2}\pi \biggl[ \sum
_{x,y,v} c_x(\eta_x) \bigl(
\nabla_{xy}\psi(\eta) \bigr)^2 \bigl(\hat{\rho
}_1\bigl(\eta, \eta^{xy}\bigr)\rho(\eta) + \hat{
\rho}_2\bigl(\eta, \eta ^{xy}\bigr)\rho\bigl(\eta
^{xy}\bigr) \bigr) \biggr]
\\
&=& - \frac{c + \delta}{4L}\pi \biggl[ \sum_{x,y}
c_x(\eta_x) \bigl(\nabla_{xy}\psi(\eta)
\bigr)^2\hat{\rho}\bigl(\eta, \eta ^{xy}\bigr) \biggr]
\\
&=& - \frac{c + \delta}{2}\mathcal{A}(\rho, \psi),
\end{eqnarray*}
so that
%
\begin{equation}
\label{eq:B2-bound} \tcB_2(\rho,\psi) \geq- \frac{c + \delta}{2}\mathcal{A}(
\rho, \psi) + I_3.
\end{equation}

Next, we observe that
\begin{eqnarray*}
\tcB_1(\rho,\psi)& = &\pi \biggl[ \sum_{\gamma, \delta}
\Gamma (\eta, \gamma , \delta) \hat{\rho}(\eta, \delta\eta)\nabla_{\delta}
\psi(\eta ) \nabla _{\gamma}\psi(\eta) \biggr]
\\
&=& \frac{1}{L^2}\pi \biggl[\sum_{x,y,v}
c_x(\eta _x) \bigl(c_x(
\eta_x) - c_x(\eta_x - 1) \bigr)
\nabla_{xy}\psi(\eta )\nabla _{xv}\psi(\eta) \hat{\rho}\bigl(
\eta, \eta^{xy}\bigr) \biggr].
\end{eqnarray*}
By another application of \eqref{eq:zrp-ass} and reversibility,
\begin{eqnarray*}
\tcB_1(\rho,\psi) &=&\frac{1}{L^2}\pi \biggl[\sum
_{x,y,v} c_x(\eta _x)
\bigl(c_x(\eta_x) - c_x(\eta_x
- 1) \bigr)\\
&&{}\times \nabla_{xy}\psi(\eta ) \bigl(\nabla_{xy}\psi(\eta)
+ \nabla_{yv}\psi\bigl(\eta^{xy}\bigr) \bigr) \hat {\rho}
\bigl(\eta , \eta^{xy}\bigr) \biggr]
\\
&\geq&2c\mathcal{A}(\rho, \psi) + \frac{1}{L^2}\pi \biggl[\sum
_
{x,y,v} c_x(\eta_x)
\bigl(c_x(\eta_x) - c_x(\eta_x
- 1) \bigr)\nabla_{xy}\psi(\eta)\\
&&{}\times \nabla_{yv}\psi\bigl(
\eta^{xy}\bigr) \hat {\rho}\bigl(\eta, \eta^{xy}\bigr) \biggr]
\\
&=:& 2c\mathcal{A}(\rho, \psi) - \frac{1}{L^2}\pi \biggl[\sum
_{x,v, y} c_x(\eta _x)
\bigl(c_y(\eta_y+1) - c_y(
\eta_y) \bigr)\\
&&{}\times \nabla_{xy}\psi(\eta )\nabla _{xv}
\psi(\eta) \hat{\rho}\bigl(\eta, \eta^{xy}\bigr) \biggr].
\end{eqnarray*}
Averaging both expressions for $\tcB_1(\rho,\psi)$, applying \eqref
{eq:zrp-ass}, and using the inequality $-2|ab| \geq-a^2 - b^2$, we obtain
%
\begin{eqnarray}
\label{eq:B1-bound} &&\tcB_1(\rho,\psi) \nonumber\\
&&\qquad\geq c \mathcal{A}(\rho, \psi)\nonumber\\
&&\qquad\quad{} +
\frac{1}{2L^2} \biggl[\sum_{x,v, y} c_x(
\eta _x) \bigl(c_x(\eta_x) -
c_x(\eta_x - 1) - \bigl(c_y(
\eta_y+1) - c_y(\eta _y) \bigr) \bigr)
\\
& &\qquad\quad{}\times\nabla_{xy}\psi(\eta)\nabla_{xv}\psi(\eta) \hat{\rho
}\bigl(\eta, \eta^{xy}\bigr) \biggr]
\nonumber\\
&&\qquad\geq \biggl(c- \frac{\delta}{2} \biggr) \mathcal{A}(\rho, \psi) -
\frac{\delta}{4L^2}\pi \biggl[\sum_{x,v, y}
c_x(\eta_x) \bigl(\nabla_{xv}\psi(\eta)
\bigr)^2\hat{\rho}\bigl(\eta, \eta^{xy}\bigr) \biggr].\nonumber
\end{eqnarray}

Adding \eqref{eq:B2-bound} and \eqref{eq:B1-bound} and applying Theorem~\ref{thmm:curv-crit}, we obtain
%
\begin{equation}
\label{eq:B-bound} \mathcal{B}(\rho, \psi) \geq \biggl(\frac{c}{2} - \delta
\biggr)\mathcal {A}(\rho, \psi) + I_4,
\end{equation}
where
\[
I_4 := I_3 - \frac{\delta}{4L^2}\pi \biggl[\sum
_{x,y,v} c_x(\eta_x) \bigl(
\nabla_{xv}\psi(\eta) \bigr)^2\hat{\rho}\bigl(\eta,
\eta^{xy}\bigr) \biggr].
\]
Since reversibility yields
\[
\pi \biggl[\sum_{x,y,v} c_x(
\eta_x) \bigl(\nabla _{xv}\psi(\eta) \bigr)^2
\hat{\rho}\bigl(\eta, \eta^{xy}\bigr) \biggr] = \pi \biggl[\sum
_{x,y,v} c_x(\eta _x) \bigl(
\nabla_{xv}\psi(\eta) \bigr)^2\hat{\rho}\bigl(
\eta^{xv}, \eta^{xy}\bigr) \biggr],
\]
we obtain after averaging,
\begin{eqnarray*}
I_4 &=& \frac{\delta}{8L^2}\pi \biggl[\sum
_{x,y,v} c_x(\eta _x) \bigl(
\nabla_{xv}\psi(\eta) \bigr)^2 \bigl( \rho\bigl(
\eta^{xy}\bigr) \bigl(\hat{\rho }_1\bigl(\eta,
\eta^{xv}\bigr) \\
&&{}+ \hat{\rho}_2\bigl(\eta,
\eta^{xv}\bigr) \bigr) - \hat {\rho }\bigl(\eta, \eta^{xy}
\bigr) - \hat{\rho}\bigl(\eta^{xv}, \eta^{xy}\bigr) \bigr)
\biggr].
\end{eqnarray*}
Applying Lemma~\ref{lem:theta_zrp} with $r = \rho(\eta^{xy})$, $s =
\rho
(\eta)$ and $t = \rho(\eta^{xv})$, we infer that
\[
I_4 \geq-\frac{\delta}{8L^2}\pi \biggl[\sum
_{x,y,v} c_x(\eta _x) \bigl(
\nabla_{xv}\psi(\eta) \bigr)^2\hat{\rho}\bigl(\eta, \eta
^{xv}\bigr) \biggr] = -\frac{\delta}{4}\mathcal{A}(\rho, \psi).
\]
In view of \eqref{eq:B-bound}, we obtain the desired result.
\end{pf*}

\subsection{Bernoulli--Laplace models}

We consider the exclusion process on the complete graph with $L$ sites
and $N$ particles, where $1 \leq N < L$. The moves are of the form
$\eta\mapsto\eta^{xy}$, where
\[
\cases{\eta^{xy}_z = \eta_z, &\quad $
\mbox{if } z \notin\{x, y \} $, \vspace *{2pt}
\cr
\eta^{xy}_x
= 0, &\quad $\mbox{if } \eta_x(1-\eta_y) = 1, \mbox{
otherwise } \eta^{xy}_x = \eta_x $,
\vspace*{2pt}
\cr
\eta^{xy}_y = 1, &\quad $\mbox{if }
\eta_x(1-\eta_y) = 1, \mbox{ otherwise }
\eta^{xy}_y = \eta_y $.} %
\]
Note that such moves conserve the total number of particles. The state
space $S$ consists of all configurations with at most one particle on
each of the $L$ sites, that is, $S = \{ \eta\in\{0,1\}^L : \sum_x
\eta_x = N \}$.
The set of moves is given by $G := \{ (x,y) \in\{1, \ldots, L\}^2 : x
\neq y \}$. To simplify notation, we will frequently write $xy$
instead of $(x,y)$. The transition rates $q : S \times G \to\R_+$ are
given by
\[
q(\eta, xy) = \frac{\lambda_x}{L}\eta_x (1-\eta_y),
\]
hence the generator $\cL$ takes the form
\[
\cL\psi(\eta) = \frac{1}{L}\sum_{xy \in G}
\lambda_x \eta_x (1-\eta _y)
\nabla_{xy} \psi(\eta),
\]
where $\nabla_{xy} \psi(\eta) = \psi(\eta^{xy}) - \psi(\eta)$.
We shall frequently use reversibility in the form
%
\begin{equation}
\label{eq:rever} \pi \biggl[ \sum_{xy \in G} F(\eta, xy) q(
\eta,xy) \biggr] = \pi \biggl[ \sum_{xy \in G} F\bigl(
\eta^{xy}, yx\bigr) q(\eta,xy) \biggr]
\end{equation}
for arbitrary functions $F : S \times G \to\R$. We observe that
\[
\mathcal{A}(\rho, \psi) = \frac{1}{2L}\pi \biggl[\sum
_{x,y} \lambda_x \bigl(\nabla_{xy}\psi(
\eta) \bigr)^2 \hat{\rho}\bigl(\eta, \eta ^{xy}\bigr)
\biggr].
\]
Note that we can omit the factor $\eta_x(1-\eta_y)$, since $\nabla_{xy}
\psi(\eta) = 0$ whenever $\eta_x(1-\eta_y) = 0$. This remark will allow
us to forget this factor in most of the expressions that we shall manipulate.

We obtain the following result on the Ricci curvature for
Bernoulli--Laplace models.

\begin{thmm}\label{thmm:BL}
Let $1 \leq N < L-1$, and assume that there exists $c > 0$ and $\delta
\in[0, 2c]$ such that
%
\begin{equation}
\label{eq:BL} c \leq\lambda_x \leq c + \delta\qquad \mbox{for all } x \in
\{ 1, \ldots, L \}.
\end{equation}
Then the Ricci curvature of the Bernoulli--Laplace model is bounded
from below by $\frac{c}{2} - \frac{7\delta}{8}$. In particular, when
$\delta< \frac{4}{7}c$, the Ricci curvature is positive.
\end{thmm}

In the special case where $\delta= 0$, we recover the result obtained
in \cite{EMT} using a different method, with exactly the same constant.

\begin{remark}\label{rem:BL-LSI}
If $\delta< \frac{4}7 c$, our result yields a modified logarithmic
Sobolev inequality with constant $c - \frac{7\delta}{4}$.
Once more, this bound is slightly weaker than the one obtained in \cite
{CDP09} (which is $c-\delta$ in our notation). Note however that both
bounds coincide for the homogeneous model, in which $\delta=0$. The
modified logarithmic Sobolev inequality for the homogeneous model has
also been studied in \cite{BT06,GQ03,Go04}.
\end{remark}

\begin{pf*}{Proof of Theorem~\ref{thmm:BL}}
As in \cite{CDP09}, we take
\[
R(\eta, xy, uv) = \cases{ %
\displaystyle \frac{1}{L^2}
\lambda_x\lambda_u\eta_x(1-
\eta_y)\eta_u(1-\eta_v), & \quad$\mbox{for }\bigl |
\{x,y,u,v\}\bigr| = 4,$
\vspace*{2pt}\cr
0, &\quad $\mbox{otherwise}.$}
\]
It follows from the definition that $\Gamma(\eta, xy, uv) = 0$ when $
|\{x,y,u,v\}| = 4$ and $\Gamma(\eta, xy, uv) = \frac{1}{L^2}\lambda
_x\lambda_u\eta_x(1-\eta_y)\eta_u(1-\eta_v)$ otherwise.
We also notice that $\nabla_{xy}\varphi\nabla_{yz} \psi= \nabla
_{xy}\varphi\nabla_{zx} \psi= 0$ for any choice of $x, y$ and $z$,
and for any $\varphi$ and $\psi$.
We then have
\begin{eqnarray*}
\tcB_2(\rho, \psi) &=& \frac{1}{2L^2}\pi \biggl[ \sum
_{x,y} \lambda_x^2 \bigl(
\nabla_{xy}\psi(\eta) \bigr)^2\hat{\rho}_1\bigl(
\eta, \eta ^{xy}\bigr)\nabla _{xy}\rho(\eta) \biggr]
\\
&&{} + \frac{1}{2L^2}\pi \biggl[\sum_{|\{x,y,u\}| = 3}
\lambda_x \lambda_u \bigl(\nabla_{xy} \psi(
\eta) \bigr)^2 \hat{\rho}_1\bigl(\eta, \eta
^{xy}\bigr)\nabla _{uy}\rho(\eta) \biggr]
\\
&&{}+ \frac{1}{2L^2}\pi \biggl[\sum_{|\{x,y,v\}| = 3}
\lambda_x^2 \bigl(\nabla _{xy} \psi(\eta)
\bigr)^2 \hat{\rho}_1\bigl(\eta, \eta^{xy}\bigr)
\nabla _{xv}\rho (\eta) \biggr]
\\
& =:& J_1+ J_2 + J_3.
\end{eqnarray*}
Using reversibility in the form of \eqref{eq:rever}, we obtain,  after
averaging with the original expression,
\begin{eqnarray*}
J_1 &=& -\frac{1}{2L^2}\pi \biggl[ \sum
_{x,y} \lambda_x \lambda_y \bigl(
\nabla _{xy}\psi(\eta) \bigr)^2\hat{\rho}_2
\bigl(\eta, \eta^{xy}\bigr)\nabla _{xy}\rho(\eta ) \biggr]
\\
&=& \frac{1}{4L^2}\pi \biggl[ \sum_{x,y}
\lambda_x \bigl(\nabla _{xy}\psi(\eta )
\bigr)^2 \bigl(\lambda_x\hat{\rho}_1\bigl(
\eta, \eta^{xy}\bigr) - \lambda_y \hat {\rho}_2
\bigl(\eta, \eta^{xy}\bigr) \bigr) \bigl(\rho\bigl(\eta^{xy}
\bigr) - \rho(\eta ) \bigr) \biggr]
\\
&\geq&-\frac{\delta}{4L^2}\pi \biggl[ \sum_{x,y}
\lambda_x \bigl(\nabla _{xy}\psi(\eta) \bigr)^2
\hat{\rho}\bigl(\eta, \eta^{xy}\bigr) \biggr] \\
&=& - \frac{\delta
}{2L}
\cA(\rho,\psi), %
\end{eqnarray*}
where the inequality is obtained using \eqref{eq:theta-cor} and \eqref{eq:BL}.

Another application of reversibility and averaging yields
\begin{eqnarray*}
J_3 &=& \frac{1}{2L^2}\pi \biggl[\sum
_{|\{x,y,v\}| = 3} \lambda_x \lambda_y \bigl(
\nabla_{xy} \psi(\eta) \bigr)^2 \hat{\rho}_2
\bigl(\eta, \eta ^{xy}\bigr) \bigl(\rho\bigl(\eta^{xv}\bigr) -
\rho\bigl(\eta^{xy}\bigr) \bigr) \biggr]
\\[-1pt]
&=& \frac{1}{4L^2}\pi \biggl[\sum_{|\{x,y,v\}| = 3}
\lambda_x \bigl(\nabla _{xy} \psi(\eta)
\bigr)^2 \bigl(\lambda_x \hat{\rho}_1\bigl(
\eta, \eta ^{xy}\bigr)\nabla_{xv}\rho(\eta)\\[-1pt]
&&{} +
\lambda_y \hat{\rho}_2\bigl(\eta, \eta ^{xy}
\bigr) \bigl(\rho\bigl(\eta^{xv}\bigr) - \rho\bigl(\eta^{xy}
\bigr) \bigr) \bigr) \biggr]
\\[-1pt]
&=& \frac{1}{4L^2}\pi \biggl[\sum_{|\{x,y,v\}| = 3}
\eta_x(1-\eta _y) (1-\eta _v)
\lambda_x \bigl(\nabla_{xy} \psi(\eta) \bigr)^2
\rho\bigl(\eta ^{xv}\bigr)\\[-1pt]
&&{}\times \bigl(\lambda_x \hat{
\rho}_1\bigl(\eta, \eta^{xy}\bigr)+ \lambda_y
\hat{\rho }_2\bigl(\eta, \eta^{xy}\bigr) \bigr) \biggr]
\\[-1pt]
&&{} - \frac{1}{4L^2}\pi \biggl[\sum_{|\{x,y,v\}| = 3}
\eta_x(1-\eta _y) (1-\eta _v)
\lambda_x \bigl(\nabla_{xy} \psi(\eta) \bigr)^2
\\[-1pt]
&&{}\times\bigl(\lambda_x \hat{\rho }_1\bigl(\eta,
\eta^{xy}\bigr)\rho(\eta)
+ \lambda_y \hat{
\rho}_2\bigl(\eta, \eta ^{xy}\bigr)\rho\bigl(
\eta^{xy}\bigr) \bigr) \biggr].
\end{eqnarray*}
Using the inequality $\lambda_x, \lambda_y \geq c \geq\delta/2$ and
Lemma~\ref{lem:theta}(i), we infer that\vspace*{-1pt}
\begin{eqnarray*}
\lambda_x \hat{\rho}_1\bigl(\eta, \eta^{xy}
\bigr) + \lambda_y \hat{\rho }_2\bigl(\eta,
\eta^{xy}\bigr)
&\geq&\frac{\delta}{2} \bigl(\hat{\rho}_1
\bigl(\eta, \eta^{xy}\bigr) + \hat {\rho }_2\bigl(\eta,
\eta^{xy}\bigr) \bigr),
\\[-1pt]
\lambda_x \hat{\rho}_1\bigl(\eta, \eta^{xy}
\bigr)\rho(\eta) + \lambda_y \hat{\rho }_2\bigl(\eta,
\eta^{xy}\bigr)\rho\bigl(\eta^{xy}\bigr)
&\geq&  c \hat{\rho}
\bigl(\eta, \eta^{xy}\bigr).
\end{eqnarray*}
Applying these bounds, we arrive at\vspace*{-1pt}
\begin{eqnarray*}
J_3 &\geq&\frac{\delta}{8L^2}\pi \biggl[\sum
_{|\{x,y,v\}| = 3} \eta _x(1-\eta _y) (1-
\eta_v)\lambda_x \bigl(\nabla_{xy} \psi(\eta)
\bigr)^2\rho \bigl(\eta ^{xv}\bigr)\\[-1pt]
&&{}\times \bigl( \hat{
\rho}_1\bigl(\eta, \eta^{xy}\bigr) + \hat{
\rho}_2\bigl(\eta, \eta ^{xy}\bigr) \bigr) \biggr]
\\[-1pt]
&&{} - \frac{(L-N -1)(c+\delta)}{2L}\mathcal{A}(\rho, \psi),
\end{eqnarray*}
where the $(L - N - 1)$ factor appears because, if there is a particle
at $x$ and no particle at $y$, there are exactly $(L - N - 1)$ possible
sites $v$ different from $x$ and $y$ where there are no
particles.

Similarly, we have
\begin{eqnarray*}
J_2 &=& \frac{1}{2L^2}\pi \biggl[\sum
_{|\{x,y,u\}| = 3} \eta_u \lambda_x
\lambda_u \bigl(\nabla_{xy} \psi(\eta) \bigr)^2
\hat{\rho}_1\bigl(\eta , \eta ^{xy}\bigr)
\nabla_{uy}\rho(\eta) \biggr]
\\
&=& \frac{1}{4L^2}\pi \biggl[\sum_{|\{x,y,u\}| = 3}
\eta_u \lambda_x \lambda _u \bigl(
\nabla_{xy} \psi(\eta) \bigr)^2\\
&&{}\times \bigl(\hat{
\rho}_1\bigl(\eta, \eta ^{xy}\bigr)\nabla_{uy}
\rho(\eta) + \hat{\rho}_2\bigl(\eta, \eta^{xy}\bigr) \bigl(
\rho \bigl(\eta^{uy}\bigr) - \rho\bigl(\eta^{xy}\bigr) \bigr)
\bigr) \biggr]
\\[-1pt]
&\geq&\frac{1}{4L^2}\pi \biggl[\sum_{|\{x,y,u\}| = 3}
\eta_u \lambda_x \lambda_u \bigl(
\nabla_{xy} \psi(\eta) \bigr)^2
\bigl(\hat{\rho
}_1\bigl(\eta, \eta^{xy}\bigr) + \hat{\rho}_2
\bigl(\eta, \eta^{xy}\bigr) \bigr)\rho\bigl(\eta^{uy}\bigr)
\biggr]
\\[-1pt]
& &{}- \frac{(N-1)(c+\delta)}{2L}\mathcal{A}(\rho, \psi)
\\[-1pt]
&\geq&-\frac{(N-1)(c+\delta)}{2L}\mathcal{A}(\rho, \psi).
\end{eqnarray*}

We now turn to $\tcB_1$. To improve readability, we shall often
suppress the variable $\eta$ in our notation. We have\vspace*{-1pt}
\begin{eqnarray*}
\tcB_1(\rho, \psi) &=& \frac{1}{L^2}\pi \biggl[\sum
_{x,y} \lambda _x^2(\nabla
_{xy}\psi)^2\hat{\rho}\bigl(\eta, \eta^{xy}\bigr)
\biggr] \\[-1pt]
&&{}+ \frac{1}{L^2}\pi \biggl[\sum_{|\{x,y,u\}| = 3}
\lambda_x \lambda_u \nabla_{xy} \psi\nabla
_{uy} \psi\hat{\rho}\bigl(\eta, \eta^{xy}\bigr) \biggr]
\\[-1pt]
& &{}+ \frac{1}{L^2}\pi \biggl[\sum_{|\{x,y,v\}| = 3}
\lambda_x^2 \nabla_{xy} \psi
\nabla_{xv}\psi\hat{\rho}\bigl(\eta, \eta^{xy}\bigr) \biggr]
\\[-1pt]
&=:& J_4 + J_5 + J_6.
\end{eqnarray*}
We have the immediate bound\vspace*{-1pt}
\[
J_4 \geq\frac{2c}{L}\mathcal{A}(\rho, \psi).
\]

Another application of \eqref{eq:rever} yields\vspace*{-1pt}
\begin{eqnarray*}
J_6 &=& \frac{1}{L^2}\pi \biggl[\sum
_{|\{x,y,v\}| = 3}\lambda_x^2\eta
_x(1-\eta_y) (1-\eta_v)
\nabla_{xy} \psi \bigl(\psi\bigl(\eta^{xv}\bigr) - \psi\bigl(
\eta ^{xy}\bigr) \bigr) \hat{\rho}\bigl(\eta, \eta^{xy}\bigr)
\biggr]
\\[-1pt]
&&{} + \frac{1}{L^2}\pi \biggl[\sum_{|\{x,y,v\}| = 3}
\lambda_x^2\eta _x(1-\eta _y)
(1-\eta_v) (\nabla_{xy} \psi )^2\hat{\rho}
\bigl(\eta, \eta ^{xy}\bigr) \biggr]
\\[-1pt]
&=&-\frac{1}{L^2}\pi \biggl[\sum_{|\{x,y,v\}| = 3}
\lambda_x\lambda _y\eta _x(1-
\eta_y) (1-\eta_v) \nabla_{xy} \psi
\nabla_{xv} \psi\hat{\rho }\bigl(\eta , \eta^{xy}\bigr)
\biggr]
\\[-1pt]
&&{} + \frac{L - N -1}{L^2}\pi \biggl[\sum_{x,y}
\lambda_x^2 \eta _x(1-\eta_y) (
\nabla_{xy} \psi )^2\hat{\rho}\bigl(\eta,
\eta^{xy}\bigr) \biggr].
\end{eqnarray*}
Averaging the latter expression with the defining formula for $J_6$, we obtain\vspace*{-1pt}
\begin{eqnarray*}
J_6 &=& \frac{L-N-1}{2L^2}\pi \biggl[\sum
_{x,y} \lambda_x^2 \eta
_x(1-\eta _y) (\nabla_{xy} \psi
)^2\hat{\rho}\bigl(\eta, \eta^{xy}\bigr) \biggr]
\\[-1pt]
&&{}+ \frac{1}{2L^2}\pi \biggl[\sum_{|\{x,y,v\}| = 3} \lambda
_x(\lambda_x -\lambda_y)
\eta_x(1-\eta_y) (1-\eta_v)
\nabla_{xy} \psi \nabla_{xv} \psi\hat{\rho}\bigl(\eta,
\eta^{xy}\bigr) \biggr].
\end{eqnarray*}
Writing
\[
J_7 := \frac{\delta}{4L^2}\pi \biggl[\sum
_{|\{x,y,v\}| = 3} \lambda_x \eta _x(1-
\eta_y) (1-\eta_v) (\nabla_{xv} \psi
)^2\hat{\rho }\bigl(\eta, \eta ^{xy}\bigr) \biggr]
\]
for brevity, we obtain using the inequality $-2|ab| \geq- a^2 - b^2$,
\begin{eqnarray*}
J_6 &\geq&\frac{c(L-N-1)}{L}\mathcal{A}(\rho,\psi) \\
&&{}-
\frac{\delta
}{4L^2}\pi \biggl[\sum_{|\{x,y,v\}| = 3}
\lambda_x \eta_x(1-\eta_y) (1-
\eta_v) (\nabla_{xy} \psi )^2\hat{\rho}\bigl(
\eta, \eta^{xy}\bigr) \biggr] - J_7
\\
&= &\frac{(2c-\delta)(L-N-1)}{2L}\mathcal{A}(\rho,\psi) - J_7.
\end{eqnarray*}
Using reversibility and averaging, we can show that
\begin{eqnarray*}
J_7 \hspace*{-1pt}&\stackrel{\scriptsize{\eqref{eq:rever}}} {=} &\hspace*{-1pt}\frac{\delta}{4L^2}\pi \biggl[
\sum_{|\{
x,y,v\}| = 3} \lambda_x
\eta_x(1-\eta_y) (1-\eta_v) (\nabla
_{xv} \psi )^2\hat{\rho}\bigl(\eta^{xv},
\eta^{xy}\bigr) \biggr]
\\
\hspace*{-2pt}&=& \frac{\delta}{8L^2}\pi \biggl[\sum_{|\{x,y,v\}| = 3}
\lambda_x \eta _x(1-\eta_y) (1-
\eta_v) (\nabla_{xv} \psi )^2 \bigl(\hat {\rho}
\bigl(\eta ^{xv}, \eta^{xy}\bigr) + \hat{\rho}\bigl(\eta,
\eta^{xy}\bigr) \bigr) \biggr]
\\
\hspace*{-4pt}&\stackrel{y \leftrightarrow v} {=}&\hspace*{-4pt} \frac{\delta}{8L^2}\pi \biggl[\sum
_{|\{
x,y,v\}| = 3} \lambda_x \eta_x(1-
\eta_y) (1-\eta_v) (\nabla _{xy} \psi
)^2 \bigl(\hat{\rho}\bigl(\eta^{xv}, \eta^{xy}
\bigr) + \hat{\rho}\bigl(\eta , \eta ^{xv}\bigr) \bigr) \biggr].
\end{eqnarray*}
Therefore, we have
\begin{eqnarray*}
J_3 + J_6 & \geq&\frac{(c-2\delta)(L-N-1)}{2L}\mathcal{A}(\rho,
\psi)
\\
&&{} + \frac{\delta}{8L^2}\pi \biggl[\sum_{|\{x,y,v\}| = 3}
\lambda_x \eta _x(1-\eta_y) (1-
\eta_v) (\nabla_{xy} \psi)^2
\\
& &{}\times \bigl(\rho\bigl(\eta^{xv}\bigr) \bigl(\hat{\rho}_1
\bigl(\eta, \eta ^{xy}\bigr) + \hat{\rho}_2\bigl(\eta,
\eta^{xy}\bigr) \bigr)-\hat{\rho}\bigl(\eta^{xv}, \eta
^{xy}\bigr) - \hat{\rho}\bigl(\eta, \eta^{xv}\bigr) \bigr)
\biggr]
\\
&\geq&\frac{(2c-5\delta)(L-N-1)}{4L}\mathcal{A}(\rho,\psi),
\end{eqnarray*}
where the last inequality has been obtained through an application of
Lemma~\ref{lem:theta_zrp}.

Similarly, we have
\begin{eqnarray*}
J_5 &\stackrel{\scriptsize{\eqref{eq:rever}}} {=}& \frac{1}{L^2}\pi \biggl[
\sum_{|\{
x,y,u\}| = 3} \lambda_x\lambda_u
\eta_x(1-\eta_y)\eta _u
\nabla_{xy} \psi\bigl(\psi(\eta) - \psi\bigl(\eta^{uy}\bigr)
\bigr) \hat{\rho }\bigl(\eta, \eta ^{xy}\bigr) \biggr]
\\
&&{}+ \frac{1}{L^2}\pi \biggl[\sum_{|\{x,y,u\}| = 3}
\lambda_x\lambda_u\eta_x(1-
\eta_y)\eta_u (\nabla_{xy} \psi
)^2\hat{\rho}\bigl(\eta, \eta^{xy}\bigr) \biggr]
\\
&\geq&- J_5 + \frac{2c(N-1)}{L}\mathcal{A}(\rho, \psi),
\end{eqnarray*}
and, therefore,
\[
J_5 \geq\frac{c(N-1)}{L}\mathcal{A}(\rho, \psi).
\]

When we sum up, we get
%
\begin{eqnarray}
\label{eq:final} \tcB(\rho, \psi) &\geq& \biggl(\frac{(L-N-1)(2c-5\delta)}{4L} +
\frac
{(N-1)(c-\delta)}{2L} + \frac{4c-\delta}{2L} \biggr)\mathcal {A}(\rho, \psi)
\nonumber
\\[-8pt]
\\[-8pt]
\nonumber
&\geq &\biggl(\frac{c}{2} - \frac{5(L-1)-3N}{4L}\delta \biggr) \mathcal {A}(
\rho, \psi).
\end{eqnarray}
At this point, what we get is a size-independent lower bound on the
Ricci curvature of $\frac{c}{2} - \frac{5\delta}{4}$. However, we can
improve this bound using the following duality argument for the
Bernoulli--Laplace process: if we consider the system with $N$
particles, and then remove all particles while simultaneously adding
particles at every empty site, we get a Bernoulli--Laplace model with
$L - N$ particles. Since this is just a change in labeling (empty and
full sites play symmetric roles), the properties of the dynamic are
invariant by this transform. Therefore, without any loss of generality,
we can assume that $N \geq L/2$. Under this extra assumption, the bound
\eqref{eq:final} leads to a lower bound on the Ricci curvature of
$\frac
{c}{2} - \frac{7\delta}{8}$, which is what we were seeking to prove.
\end{pf*}

\subsection{The random transposition model}

We now consider the random transposition model, which is a random walk
on the group of permutations $\mathcal{S}_n$ with $n \geq2$. If we
denote by $\mathcal{T}_n$ the set of all transpositions in $\mathcal
{S}_n$, we can write the generator of the dynamics as
\[
\mathcal{L}f(\sigma) := \frac{2}{n(n-1)} \sum_{\tau\in\mathcal
{T}_n}
\nabla_{\tau}f(\sigma),
\]
where $\nabla_{\tau} f(\sigma) = f(\tau\circ\sigma) - f(\sigma)$ and
$\pi$ is the uniform measure on $\mathcal{S}_n$, that is, $\pi
(\sigma)
= (n!)^{-1}$ for all $\sigma\in\mathcal{S}_n$. We write $\tau=
(i,j)$ for the transposition that swaps $i$ and $j$. We also write
$(i,j,k)$ for the mapping that cyclically permutes $i,j$ and $k$. To
simplify notation, we will use the shorthand notation
\[
\sigma_{ij} := (i,j) \circ\sigma,\qquad \sigma_{ijk} := (i,j,k)
\circ\sigma, \qquad \nabla_{ij} := \nabla_{(i,j)}.
\]
We obtain the following result.
%
\begin{thmm}\label{th4.9}
For $n \geq2$, the Ricci curvature of the random transposition model
is bounded from below by $\frac{4}{n(n-1)}$.
\end{thmm}

This bound was already obtained using a different method in \cite{EMT}.
It should be noted that the modified logarithmic Sobolev constant
implied by this Ricci curvature bound is significantly worse (by a
factor $1/n$) than its known optimal behaviour. We do not know whether
the Ricci curvature's true behaviour should match the MLSI-constant.

\begin{pf*}{Proof of Theorem \ref{th4.9}}
It follows from the definition that the action functional can be
written as
\[
\mathcal{A}(\rho, \psi) = \frac{1}{2n(n-1)}\pi \biggl[\sum
_{i \neq
j} \bigl(\nabla_{ij}\psi(\sigma)
\bigr)^2\hat{\rho}(\sigma, \sigma _{ij}) \biggr].
\]
Note that a factor $\frac{1}2$ appears, since every transposition $(i,j)$
is counted twice. We define $R$ as
\[
R\bigl(\sigma, (i,j), (k,\ell)\bigr) = \cases{ \displaystyle
\frac{4}{n^2(n-1)^2}, &\quad $\mbox{if } \bigl|\{i,j,k,\ell\}\bigr| = 4,$
\vspace*{2pt}\cr
0, & \quad$\mbox{otherwise}.$}
\]
Using reversibility and the fact that $(i,j)^{-1} = (i,j)$, it follows that
\begin{eqnarray*}
\tcB_2(\rho, \psi) &=& \frac{2}{n^2(n-1)^2}\pi \biggl[\sum
_{|\{
i,j,k\}
|=3} \bigl(\nabla_{ij}\psi(\sigma)
\bigr)^2\hat{\rho}_1(\sigma, \sigma _{ij})
\nabla_{ik}\rho(\sigma) \biggr]
\\
&&{} + \frac{2}{n^2(n-1)^2}\pi \biggl[\sum_{i\neq j} \bigl(
\nabla _{ij}\psi (\sigma) \bigr)^2\hat{\rho}_1(
\sigma, \sigma_{ij})\nabla_{ij}\rho (\sigma ) \biggr]
\\
&=& \frac{2}{n^2(n-1)^2}\pi \biggl[\sum_{|\{i,j,k\}|=3} \bigl(
\nabla _{ij}\psi(\sigma) \bigr)^2\hat{\rho}_2(
\sigma, \sigma_{ij}) \bigl(\rho (\sigma_{ijk}) - \rho(
\sigma_{ij}) \bigr) \biggr]
\\
&&{} - \frac{1}{n^2(n-1)^2}\pi \biggl[\sum_{i\neq j} \bigl(
\nabla _{ij}\psi (\sigma) \bigr)^2\hat{\rho}_2(
\sigma, \sigma_{ij})\nabla_{ij}\rho (\sigma ) \biggr].
\end{eqnarray*}
Averaging the latter two expressions, we obtain
\begin{eqnarray*}
\tcB_2(\rho, \psi)
&=& \frac{1}{n^2(n-1)^2}\pi \biggl[\sum
_{|\{
i,j,k\}|=3} \bigl(\nabla_{ij}\psi(\sigma)
\bigr)^2\\
&&{}\times \bigl(\rho(\sigma_{ik})\hat{\rho}_1(
\sigma, \sigma_{ij})+ \rho (\sigma _{ijk})\hat{
\rho}_2(\sigma, \sigma_{ij}) \bigr) \biggr]
\\
& &{}- \frac{1}{n^2(n-1)^2}\pi \biggl[\sum_{|\{i,j,k\}|=3} \bigl(
\nabla_{ij}\psi(\sigma) \bigr)^2 \\
&&{}\times\bigl(\rho(\sigma)\hat{
\rho}_1(\sigma, \sigma_{ij}) + \rho (\sigma _{ij})
\hat{\rho}_2(\sigma, \sigma_{ij}) \bigr) \biggr]
\\
&&{} +\frac{1}{2 n^2(n-1)^2}\pi \biggl[\sum_{i\neq j} \bigl(
\nabla _{ij}\psi (\sigma) \bigr)^2 \bigl(\hat{
\rho}_1(\sigma, \sigma_{ij})-\hat {\rho }_2(
\sigma, \sigma_{ij}) \bigr)\nabla_{ij}\rho(\sigma) \biggr].
\end{eqnarray*}
Using (i) and (ii) of Lemma~\ref{lem:theta} and \eqref
{eq:theta-cor}, we infer that
\[
\tcB_2(\rho, \psi) \geq\frac{1}{n^2(n-1)^2}\pi \biggl[\sum
_{|\{i,j,k\}|=3} \bigl(\nabla _{ij}\psi(\sigma)
\bigr)^2\hat{\rho}(\sigma_{ik}, \sigma _{ijk})
\biggr] - \frac{2(n-2)}{n(n-1)}\mathcal{A}(\rho, \psi).
\]
The term $\tcB_1(\rho, \psi)$ can be written as
\begin{eqnarray*}
\tcB_1(\rho, \psi) &=& \frac{2}{n^2(n-1)^2}\pi \biggl[ \sum
_{i\neq j} \bigl(\nabla_{ij}\psi(\sigma)
\bigr)^2\hat{\rho}(\sigma, \sigma _{ij}) \biggr]
\\
&&{} + \frac{4}{n^2(n-1)^2}\pi \biggl[\sum_{|\{i,j,k\}|=3}
\nabla_{ij}\psi(\sigma)\nabla_{ik}\psi(\sigma)\hat{\rho}(
\sigma , \sigma _{ij}) \biggr]
\\
&=&\frac{4}{n(n-1)}\mathcal{A}(\rho, \psi) + \frac{4}{n^2(n-1)^2} \pi \biggl[
\sum_{|\{i,j,k\}|=3} \bigl(\nabla_{ij}\psi(\sigma)
\bigr)^2\hat{\rho}(\sigma, \sigma _{ij}) \biggr]
\\
&&{} + \frac{4}{n^2(n-1)^2}\pi \biggl[ \sum_{|\{i,j,k\}|=3} \nabla
_{ij}\psi (\sigma) \bigl(\psi(\sigma_{ik}) - \psi(
\sigma_{ij}) \bigr)\hat {\rho}(\sigma , \sigma_{ij}) \biggr].
\end{eqnarray*}
Using reversibility and averaging, the latter term can be reformulated as
\begin{eqnarray*}
&&\frac{4}{n^2(n-1)^2}\pi \biggl[ \sum_{|\{i,j,k\}|=3}
\nabla_{ij}\psi (\sigma) \bigl(\psi(\sigma_{ik}) - \psi(
\sigma_{ij}) \bigr)\hat {\rho}(\sigma , \sigma_{ij}) \biggr]
\\
&&\qquad = \frac{4}{n^2(n-1)^2}\pi \biggl[ \sum_{|\{i,j,k\}|=3} \nabla
_{ij}\psi (\sigma) \bigl(\psi(\sigma) - \psi(\sigma_{ijk})
\bigr)\hat{\rho }(\sigma, \sigma_{ij}) \biggr]
\\
&&\qquad = \frac{2}{n^2(n-1)^2}\pi \biggl[ \sum_{|\{i,j,k\}|=3} \nabla
_{ij}\psi (\sigma) \bigl(\psi(\sigma_{ik}) - \psi(
\sigma_{ij}) + \psi(\sigma ) - \psi (\sigma_{ijk}) \bigr)\\
&&\qquad\quad{}\times\hat{
\rho}(\sigma, \sigma_{ij}) \biggr]
\\
&&\qquad = - \frac{4(n-2)}{n(n-1)}\cA(\rho, \psi)\\
&&\qquad\quad{} + \frac{2}{n^2(n-1)^2}\pi \biggl[ \sum
_{|\{i,j,k\}|=3} \nabla _{ij}\psi (\sigma) \bigl(
\psi(\sigma_{ik}) - \psi(\sigma_{ijk}) \bigr)\hat {\rho }(
\sigma, \sigma_{ij}) \biggr],
\end{eqnarray*}
which yields
\begin{eqnarray*}
\tcB_1(\rho, \psi) 
&=&
\biggl(\frac{4}{n(n-1)} + \frac{4(n-2)}{n(n-1)} \biggr)\mathcal {A}(\rho , \psi)
\\
&&{} + \frac{2}{n^2(n-1)^2}\pi \biggl[ \sum_{|\{i,j,k\}|=3} \nabla
_{ij}\psi(\sigma) \bigl(\psi(\sigma_{ik}) - \psi(
\sigma_{ijk}) \bigr)\hat {\rho}(\sigma, \sigma_{ij}) \biggr].
\end{eqnarray*}
Young's inequality and reversibility imply that
\begin{eqnarray*}
\tcB_1(\rho, \psi)&\geq& \biggl(\frac{4}{n} -
\frac
{2(n-2)}{n(n-1)} \biggr)\mathcal{A}(\rho, \psi) \\
&&{}- \frac{1}{n^2(n-1)^2}\pi \biggl[
\sum_{|\{
i,j,k\}
|=3} \bigl(\psi(\sigma_{ik}) - \psi(
\sigma_{ijk}) \bigr)^2\hat{\rho }(\sigma ,
\sigma_{ij}) \biggr]
\\
&=&\frac{2}{n-1}\mathcal{A}(\rho, \psi) - \frac{1}{n^2(n-1)^2}\pi \biggl[ \sum
_{|\{i,j,k\}|=3} \bigl(\nabla_{jk}\psi(\sigma)
\bigr)^2\hat {\rho }(\sigma_{ik}, \sigma_{ikj})
\biggr]
\\
&\mathop{=}\limits_{k \leftrightarrow i}&\frac{2}{n-1}\mathcal{A}(\rho, \psi) -
\frac{1}{n^2(n-1)^2}\pi \biggl[ \sum_{|\{i,j,k\}|=3} \bigl(\nabla
_{ij}\psi (\sigma) \bigr)^2\hat{\rho}(\sigma_{ik},
\sigma_{ijk}) \biggr].
\end{eqnarray*}
When we sum $\tcB_1$ and $\tcB_2$, we get
\[
\mathcal{B}(\rho, \psi) \geq\frac{4}{n(n-1)}\mathcal{A}(\rho, \psi),
\]
which is the desired result.
\end{pf*}

\begin{appendix}\label{app}
\section*{Appendix: Properties of the logarithmic mean}
In this paper, we make use of some basic properties of the logarithmic
mean, given for $a,b \geq0$ (resp., $a,b > 0$) by
\[
\theta(a,b) := \int_0^1{a^{1-p}b^p
\,\mathrm{ d} p} = \frac{a - b}{\log a -
\log b}.
\]
The following properties are taken from \cite{EM12}; see Lemma~5.4. We
write $\theta_i(s,t) = \partial_i \theta(s,t)$ for $i=1,2$.
%
\begin{lemm}\label{lem:theta}
The following assertions hold:
\begin{longlist}[(iii)]
\item[(i)] $u\theta_1(s,t) + v\theta_2(s,t) \geq\theta(u,v) $ for
all $s, t, u, v > 0$;

\item[(ii)] $s\theta_1(s,t) + t\theta_2(s,t) = \theta(s,t)$ for all
$s, t > 0$;

\item[(iii)] $\theta$ is symmetric, concave and increasing in both variables;

\item[(iv)] $\theta$ is positively 1-homogeneous, that is, for all
$\lambda, s, t \geq0$, we have $\theta(\lambda s, \lambda t) =
\lambda
\theta(s, t)$.
\end{longlist}
\end{lemm}

It follows directly from this lemma that for every $\lambda_1, \lambda
_2 > 0$ and $s, t \geq0$,
%
\begin{equation}
\label{eq:theta-cor} \bigl(\lambda_1 \theta_1(s,t)
- \lambda_2\theta_2(s,t)\bigr) (s-t) \leq \bigl( \max\{
\lambda_1, \lambda_2 \} - \min\{\lambda_1,
\lambda_2 \} \bigr) \theta(s,t). %
\end{equation}
The following inequality plays a crucial role in the proof of the
curvature bound for the Bernoulli--Laplace model and for zero-range processes.

\begin{lemm}\label{lem:theta_zrp}
For any $r \geq0$ and $s, t \geq0$, we have
\[
r \bigl(\theta_1(s,t) + \theta_2(s,t) \bigr) - \bigl(
\theta(r,s) + \theta (r,t) \bigr) \geq- \theta(s,t).
\]
\end{lemm}

\begin{pf}
If $r = 0$, the inequality is trivially true, so without loss of
generality we can assume that $r > 0$. Let $u = s/r$ and $v = t/r$.
Using the fact that $\theta$ is $1$-homogeneous, and $\theta_1$ and
$\theta_2$ are $0$-homogeneous, the inequality we wish to prove is
equivalent to
%
\begin{equation}
\label{eq:lem_concavity_zrp} \theta_1(u,v) + \theta_2(u,v) -
\theta(1,u) - \theta(1,v) \geq -\theta (u,v).
\end{equation}
Since $\theta$ is concave, we have
\begin{eqnarray*}
\theta(1,u) + \theta(1,v) &=& 2 \times\frac{1}{2} \bigl(\theta(u,1) +
\theta(1,v) \bigr)
\\
&\leq&2\theta \biggl(\frac{u+1}{2}, \frac{v+1}{2} \biggr)
\\
&= &\theta(u+1, v+1).
\end{eqnarray*}
Using the ``curve below tangent'' formulation of concavity applied to the
function $x \mapsto\theta(u + x, v+x)$, we have
\[
\theta(u+1, v+1) \leq\theta(u, v) + \theta_1(u,v) +
\theta_2(u,v),
\]
and \eqref{eq:lem_concavity_zrp} immediately follows, which completes
the proof.
\end{pf}
\end{appendix}

\section*{Acknowledgements} Part of this work has been done
while M. Fathi visited J.~Maas at the University of Bonn in July 2014. We would
like to thank the referees for their careful reading of the manuscript.

%

\printaddresses

\begin{thebibliography}{24}

\bibitem{BE85}
%
\begin{bincollection}[mr]
\bauthor{\bsnm{Bakry},~\bfnm{D.}\binits{D.}} \AND
\bauthor{\bsnm{{\'E}mery},~\bfnm{Michel}\binits{M.}}
(\byear{1985}).
\btitle{Diffusions hypercontractives}.
In \bbooktitle{S\'eminaire de Probabilit\'es, XIX, 1983/84}.
\bseries{Lecture Notes in Math.}
\bvolume{1123}
\bpages{177--206}.
\bpublisher{Springer},
\blocation{Berlin}.
\bid{doi={10.1007/BFb0075847}, mr={0889476}}
\end{bincollection}
%
\bptok{imsref}%
\endbibitem

\bibitem{BT06}
%
\begin{barticle}[mr]
\bauthor{\bsnm{Bobkov},~\bfnm{Sergey~G.}\binits{S.~G.}} \AND
\bauthor{\bsnm{Tetali},~\bfnm{Prasad}\binits{P.}}
(\byear{2006}).
\btitle{Modified logarithmic {S}obolev inequalities in discrete settings}.
\bjournal{J. Theoret. Probab.}
\bvolume{19}
\bpages{289--336}.
\bid{doi={10.1007/s10959-006-0016-3}, issn={0894-9840}, mr={2283379}}
\end{barticle}
%
\bptok{imsref}%
\endbibitem

\bibitem{BS09}
%
\begin{barticle}[mr]
\bauthor{\bsnm{Bonciocat},~\bfnm{Anca-Iuliana}\binits{A.-I.}} \AND
\bauthor{\bsnm{Sturm},~\bfnm{Karl-Theodor}\binits{K.-T.}}
(\byear{2009}).
\btitle{Mass transportation and rough curvature bounds for discrete spaces}.
\bjournal{J. Funct. Anal.}
\bvolume{256}
\bpages{2944--2966}.
\bid{doi={10.1016/j.jfa.2009.01.029}, issn={0022-1236}, mr={2502429}}
\end{barticle}
%
\bptok{imsref}%
\endbibitem

\bibitem{BCDP06}
%
\begin{barticle}[mr]
\bauthor{\bsnm{Boudou},~\bfnm{Anne-Severine}\binits{A.-S.}},
\bauthor{\bsnm{Caputo},~\bfnm{Pietro}\binits{P.}},
\bauthor{\bsnm{Dai Pra},~\bfnm{Paolo}\binits{P.}} \AND
\bauthor{\bsnm{Posta},~\bfnm{Gustavo}\binits{G.}}
(\byear{2006}).
\btitle{Spectral gap estimates for interacting particle systems via a
{B}ochner-type identity}.
\bjournal{J. Funct. Anal.}
\bvolume{232}
\bpages{222--258}.
\bid{doi={10.1016/j.jfa.2005.07.012}, issn={0022-1236}, mr={2200172}}
\bptnote{check year}%
\end{barticle}
%
\bptok{imsref}%
\endbibitem

\bibitem{CDP09}
%
\begin{barticle}[mr]
\bauthor{\bsnm{Caputo},~\bfnm{Pietro}\binits{P.}},
\bauthor{\bsnm{Dai Pra},~\bfnm{Paolo}\binits{P.}} \AND
\bauthor{\bsnm{Posta},~\bfnm{Gustavo}\binits{G.}}
(\byear{2009}).
\btitle{Convex entropy decay via the {B}ochner--{B}akry--{E}mery approach}.
\bjournal{Ann. Inst. Henri Poincar\'e Probab. Stat.}
\bvolume{45}
\bpages{734--753}.
\bid{doi={10.1214/08-AIHP183}, issn={0246-0203}, mr={2548501}}
\bptnote{check pages}%
\end{barticle}
%
\bptok{imsref}%
\endbibitem

\bibitem{CP07}
%
\begin{barticle}[mr]
\bauthor{\bsnm{Caputo},~\bfnm{Pietro}\binits{P.}} \AND
\bauthor{\bsnm{Posta},~\bfnm{Gustavo}\binits{G.}}
(\byear{2007}).
\btitle{Entropy dissipation estimates in a zero-range dynamics}.
\bjournal{Probab. Theory Related Fields}
\bvolume{139}
\bpages{65--87}.
\bid{doi={10.1007/s00440-006-0039-9}, issn={0178-8051}, mr={2322692}}
\end{barticle}
%
\bptok{imsref}%
\endbibitem

\bibitem{CHLZ}
%
\begin{barticle}[mr]
\bauthor{\bsnm{Chow},~\bfnm{Shui-Nee}\binits{S.-N.}},
\bauthor{\bsnm{Huang},~\bfnm{Wen}\binits{W.}},
\bauthor{\bsnm{Li},~\bfnm{Yao}\binits{Y.}} \AND
\bauthor{\bsnm{Zhou},~\bfnm{Haomin}\binits{H.}}
(\byear{2012}).
\btitle{Fokker--{P}lanck equations for a free energy functional or
{M}arkov process on a graph}.
\bjournal{Arch. Ration. Mech. Anal.}
\bvolume{203}
\bpages{969--1008}.
\bid{issn={0003-9527}, mr={2928139}}
\end{barticle}
%
\bptok{imsref}%
\endbibitem

\bibitem{DP12}
%
\begin{barticle}[mr]
\bauthor{\bsnm{Dai Pra},~\bfnm{Paolo}\binits{P.}} \AND
\bauthor{\bsnm{Posta},~\bfnm{Gustavo}\binits{G.}}
(\byear{2013}).
\btitle{Entropy decay for interacting systems via the
{B}ochner--{B}akry--\'Emery approach}.
\bjournal{Electron. J. Probab.}
\bvolume{18}
\bpages{52}.
\bid{doi={10.1214/EJP.v18-2041}, issn={1083-6489}, mr={3065862}}
\bptnote{check volume, check pages, check year}%
\end{barticle}
%
\bptok{imsref}%
\endbibitem

\bibitem{EM12}
%
\begin{barticle}[mr]
\bauthor{\bsnm{Erbar},~\bfnm{Matthias}\binits{M.}} \AND
\bauthor{\bsnm{Maas},~\bfnm{Jan}\binits{J.}}
(\byear{2012}).
\btitle{Ricci curvature of finite {M}arkov chains via convexity of the
entropy}.
\bjournal{Arch. Ration. Mech. Anal.}
\bvolume{206}
\bpages{997--1038}.
\bid{doi={10.1007/s00205-012-0554-z}, issn={0003-9527}, mr={2989449}}
\end{barticle}
%
\bptok{imsref}%
\endbibitem

\bibitem{EM13}
%
\begin{barticle}[mr]
\bauthor{\bsnm{Erbar},~\bfnm{Matthias}\binits{M.}} \AND
\bauthor{\bsnm{Maas},~\bfnm{Jan}\binits{J.}}
(\byear{2014}).
\btitle{Gradient flow structures for discrete porous medium equations}.
\bjournal{Discrete Contin. Dyn. Syst.}
\bvolume{34}
\bpages{1355--1374}.
\bid{issn={1078-0947}, mr={3117845}}
\end{barticle}
%
\bptok{imsref}%
\endbibitem

\bibitem{EMT}
%
\begin{barticle}[mr]
\bauthor{\bsnm{Erbar},~\bfnm{Matthias}\binits{M.}},
\bauthor{\bsnm{Maas},~\bfnm{Jan}\binits{J.}} \AND
\bauthor{\bsnm{Tetali},~\bfnm{P.}\binits{P.}}
(\byear{2015}).
\btitle{Discrete Ricci curvature bounds for Bernoulli--Laplace and random
transposition models}.
\bjournal{Annales Fac. Sci. Toulouse (6)}
\bvolume{24}
\bpages{781--800}.
\bid{mr={3434256}}
\end{barticle}
%
\bptok{imsref}%
\endbibitem

\bibitem{GQ03}
%
\begin{barticle}[mr]
\bauthor{\bsnm{Gao},~\bfnm{Fuqing}\binits{F.}} \AND
\bauthor{\bsnm{Quastel},~\bfnm{Jeremy}\binits{J.}}
(\byear{2003}).
\btitle{Exponential decay of entropy in the random transposition and
{B}ernoulli--{L}aplace models}.
\bjournal{Ann. Appl. Probab.}
\bvolume{13}
\bpages{1591--1600}.
\bid{doi={10.1214/aoap/1069786512}, issn={1050-5164}, mr={2023890}}
\end{barticle}
%
\bptok{imsref}%
\endbibitem

\bibitem{GM12}
%
\begin{barticle}[mr]
\bauthor{\bsnm{Gigli},~\bfnm{Nicola}\binits{N.}} \AND
\bauthor{\bsnm{Maas},~\bfnm{Jan}\binits{J.}}
(\byear{2013}).
\btitle{Gromov--{H}ausdorff convergence of discrete transportation metrics}.
\bjournal{SIAM J. Math. Anal.}
\bvolume{45}
\bpages{879--899}.
\bid{doi={10.1137/120886315}, issn={0036-1410}, mr={3045651}}
\end{barticle}
%
\bptok{imsref}%
\endbibitem

\bibitem{Go04}
%
\begin{barticle}[mr]
\bauthor{\bsnm{Goel},~\bfnm{Sharad}\binits{S.}}
(\byear{2004}).
\btitle{Modified logarithmic {S}obolev inequalities for some models of
random walk}.
\bjournal{Stochastic Process. Appl.}
\bvolume{114}
\bpages{51--79}.
\bid{doi={10.1016/j.spa.2004.06.001}, issn={0304-4149}, mr={2094147}}
\end{barticle}
%
\bptok{imsref}%
\endbibitem

\bibitem{JKO}
%
\begin{barticle}[mr]
\bauthor{\bsnm{Jordan},~\bfnm{Richard}\binits{R.}},
\bauthor{\bsnm{Kinderlehrer},~\bfnm{David}\binits{D.}} \AND
\bauthor{\bsnm{Otto},~\bfnm{Felix}\binits{F.}}
(\byear{1998}).
\btitle{The variational formulation of the {F}okker--{P}lanck equation}.
\bjournal{SIAM J. Math. Anal.}
\bvolume{29}
\bpages{1--17}.
\bid{doi={10.1137/S0036141096303359}, issn={0036-1410}, mr={1617171}}
\end{barticle}
%
\bptok{imsref}%
\endbibitem

\bibitem{LY10}
%
\begin{barticle}[mr]
\bauthor{\bsnm{Lin},~\bfnm{Yong}\binits{Y.}} \AND
\bauthor{\bsnm{Yau},~\bfnm{Shing-Tung}\binits{S.-T.}}
(\byear{2010}).
\btitle{Ricci curvature and eigenvalue estimate on locally finite graphs}.
\bjournal{Math. Res. Lett.}
\bvolume{17}
\bpages{343--356}.
\bid{doi={10.4310/MRL.2010.v17.n2.a13}, issn={1073-2780}, mr={2644381}}
\end{barticle}
%
\bptok{imsref}%
\endbibitem

\bibitem{LV09}
%
\begin{barticle}[mr]
\bauthor{\bsnm{Lott},~\bfnm{John}\binits{J.}} \AND
\bauthor{\bsnm{Villani},~\bfnm{C{\'e}dric}\binits{C.}}
(\byear{2009}).
\btitle{Ricci curvature for metric-measure spaces via optimal transport}.
\bjournal{Ann. of Math. (2)}
\bvolume{169}
\bpages{903--991}.
\bid{doi={10.4007/annals.2009.169.903}, issn={0003-486X}, mr={2480619}}
\end{barticle}
%
\bptok{imsref}%
\endbibitem

\bibitem{Maas11}
%
\begin{barticle}[mr]
\bauthor{\bsnm{Maas},~\bfnm{Jan}\binits{J.}}
(\byear{2011}).
\btitle{Gradient flows of the entropy for finite {M}arkov chains}.
\bjournal{J. Funct. Anal.}
\bvolume{261}
\bpages{2250--2292}.
\bid{doi={10.1016/j.jfa.2011.06.009}, issn={0022-1236}, mr={2824578}}
\end{barticle}
%
\bptok{imsref}%
\endbibitem

\bibitem{Mie11}
%
\begin{barticle}[mr]
\bauthor{\bsnm{Mielke},~\bfnm{Alexander}\binits{A.}}
(\byear{2011}).
\btitle{A gradient structure for reaction--diffusion systems and for
energy-drift-diffusion systems}.
\bjournal{Nonlinearity}
\bvolume{24}
\bpages{1329--1346}.
\bid{doi={10.1088/0951-7715/24/4/016}, issn={0951-7715}, mr={2776123}}
\end{barticle}
%
\bptok{imsref}%
\endbibitem

\bibitem{Mi13}
%
\begin{barticle}[mr]
\bauthor{\bsnm{Mielke},~\bfnm{Alexander}\binits{A.}}
(\byear{2013}).
\btitle{Geodesic convexity of the relative entropy in reversible
{M}arkov chains}.
\bjournal{Calc. Var. Partial Differential Equations}
\bvolume{48}
\bpages{1--31}.
\bid{doi={10.1007/s00526-012-0538-8}, issn={0944-2669}, mr={3090532}}
\end{barticle}
%
\bptok{imsref}%
\endbibitem

\bibitem{Mor}
%
\begin{barticle}[mr]
\bauthor{\bsnm{Morris},~\bfnm{Ben}\binits{B.}}
(\byear{2006}).
\btitle{Spectral gap for the zero range process with constant rate}.
\bjournal{Ann. Probab.}
\bvolume{34}
\bpages{1645--1664}.
\bid{doi={10.1214/009117906000000304}, issn={0091-1798}, mr={2271475}}
\end{barticle}
%
\bptok{imsref}%
\endbibitem

\bibitem{Oll09}
%
\begin{barticle}[mr]
\bauthor{\bsnm{Ollivier},~\bfnm{Yann}\binits{Y.}}
(\byear{2009}).
\btitle{Ricci curvature of {M}arkov chains on metric spaces}.
\bjournal{J. Funct. Anal.}
\bvolume{256}
\bpages{810--864}.
\bid{doi={10.1016/j.jfa.2008.11.001}, issn={0022-1236}, mr={2484937}}
\end{barticle}
%
\bptok{imsref}%
\endbibitem

\bibitem{Oll12}
%
\begin{bincollection}[mr]
\bauthor{\bsnm{Ollivier},~\bfnm{Yann}\binits{Y.}}
(\byear{2013}).
\btitle{A visual introduction to {R}iemannian curvatures and some
discrete generalizations}.
In \bbooktitle{Analysis and Geometry of Metric Measure Spaces}.
\bseries{CRM Proc. Lecture Notes}
\bvolume{56}
\bpages{197--220}.
\bpublisher{Amer. Math. Soc.},
\blocation{Providence, RI}.
\bid{mr={3060504}}
\bptnote{check pages, check year}%
\end{bincollection}
%
\bptok{imsref}%
\endbibitem

\bibitem{St06}
%
\begin{barticle}[mr]
\bauthor{\bsnm{Sturm},~\bfnm{Karl-Theodor}\binits{K.-T.}}
(\byear{2006}).
\btitle{On the geometry of metric measure spaces. {I}}.
\bjournal{Acta Math.}
\bvolume{196}
\bpages{65--131}.
\bid{doi={10.1007/s11511-006-0002-8}, issn={0001-5962}, mr={2237206}}
\end{barticle}
%
\bptok{imsref}%
\endbibitem
\end{thebibliography}
\end{document}